\documentclass[12pt]{amsart}
\usepackage[alphabetic]{amsrefs}

\usepackage{cite}
\usepackage{hyperref}

\usepackage[margin=0.75in]{geometry}
\usepackage{graphicx,epstopdf,color}
\newtheorem{theorem}{Theorem}[section]
\newtheorem{lemma}[theorem]{Lemma}
\numberwithin{equation}{section}
\newcommand{\R}{\mathbb{R}}
\newcommand{\Per}{\mbox{\rm{Per}}}

\usepackage{amssymb}
\renewcommand{\leq}{\leqslant}
\renewcommand{\le}{\leqslant}

\renewcommand{\ge}{\geqslant}

 \author{Benjamin Baronowitz}
 \author{Serena Dipierro}
 \author{Enrico Valdinoci}
 
 \address{Department of Mathematics and Statistics
 \newline\indent University of Western Australia \newline\indent
 35 Stirling Highway, WA 6009 Crawley, Australia}
 \email{22500237@student.uwa.edu.au, serena.dipierro@uwa.edu.au, enrico.valdinoci@uwa.edu.au}
 
\makeatletter
\@namedef{subjclassname@2020}{\textup{2020} Mathematics Subject Classification}
\makeatother
\subjclass[2020]{35R11, 49Q05}

\keywords{Nonlocal minimal surfaces, stickiness, qualitative and quantitative behavior.}

\thanks{S. Dipierro and E. Valdinoci are members of AustMS.
S. Dipierro is supported by
the Australian Research Council DECRA DE180100957
``PDEs, free boundaries and applications''.
E. Valdinoci is supported by the Australian Laureate Fellowship
FL190100081
``Minimal surfaces, free boundaries and partial differential equations''.}

\title[The stickiness property for antisymmetric nonlocal minimal graphs]{The stickiness property \\ for antisymmetric nonlocal minimal graphs}

\begin{document}

\begin{abstract}
We show that arbitrarily small antisymmetric perturbations of the zero function are sufficient to
produce the stickiness phenomenon for planar nonlocal minimal graphs
(with the same quantitative bounds obtained for the case of even symmetric perturbations, up to multiplicative constants).

In proving this result, one also establishes an odd symmetric version of the maximum principle
for nonlocal minimal graphs, according to which the odd symmetric minimizer is positive in the direction of the positive bump and negative in the direction of the negative bump.
\end{abstract}

\maketitle

\section{Introduction}

Nonlocal minimal surfaces were introduced in~\cite{MR2675483}
as minimizers of a fractional perimeter functional with respect to some given external datum.
As established in~\cite{MR3516886}, when the external datum is a graph in some direction, so is the whole minimizer,
hence it is also common to consider the case of nonlocal minimal graphs, i.e. of nonlocal minimal surfaces which
possess a graphical structure.

Nonlocal minimal graphs exhibit a quite peculiar phenomenon, discovered
in~\cite{MR3596708}, called ``stickiness''.
Roughly speaking, different from the classical minimal surfaces, in the fractional setting
remote interactions are capable of producing boundary discontinuities, which in turn
are complemented by boundary divergence of the derivative of the nonlocal minimal
graph: more specifically, boundary discontinuities of nonlocal
minimal graphs are equivalent to the divergence at the boundary of the first derivative
of the graph, see~\cite[Corollary~1.3]{MR4104542},
and this can also be interpreted as a ``butterfly effect'', since a small perturbation 
of the boundary datum not only produces a small discontinuity of the graph at the
boundary but it also forces
the slope of the graph at the boundary to shift at once from a finite
value (even zero) to infinity, see~\cite[Figure~3]{MR4096831}.

We focus here on the case in which the ambient space is of dimension~$2$: this is indeed
the simplest possible stage to detect interesting geometric patterns and, differently from
the classical case, it already provides a number of difficulties
since in the nonlocal framework segments are in general not the boundary of nonlocal
minimal objects.
Also, in the plane the stickiness phenomenon is known to be essentially ``generic''
with respect to the external datum, see~\cite{MR4104542}, in the sense that boundary
discontinuities and boundary singularities for the derivatives
can be produced by arbitrarily small perturbations of a given external datum. In particular,
these behaviors can be obtained by arbitrarily small perturbations of the flat case in which
the external datum is, say,
the subgraph of the function which is identically zero.

See also~\cite{MR4178752} for an analysis of nonlocal minimal graphs
in dimension~$3$, and~\cite{MR3926519, MR4184583,
2020arXiv201000798D}
for other examples of stickiness.

Up to now, all the examples of stickiness in the setting
of graphs were constructed by using ``one side bumps'' in the perturbation
(for instance, adding suitable positive bumps to the zero function, or to a given function
outside a vertical slab).
In this paper, we construct examples of stickiness in which the datum is antisymmetric.
In a nutshell, we will consider arbitrarily small perturbations of the zero function by
bumps possessing odd symmetry: in this case, in principle, one may fear that the effects of
equal and opposite bumps would cancel each other and prevent the stickiness phenomenon
to occur, but we will instead establish that the stickiness
phenomenon is persistent also in this class of antisymmetric perturbations of the flat case.

Also, we provide quantitative bounds on the resulting boundary discontinuity which turn out 
to be as good as the ones
available for positive bumps (up to multiplicative constants).\medskip

\begin{figure}
		\centering
		\includegraphics[width=.4\linewidth]{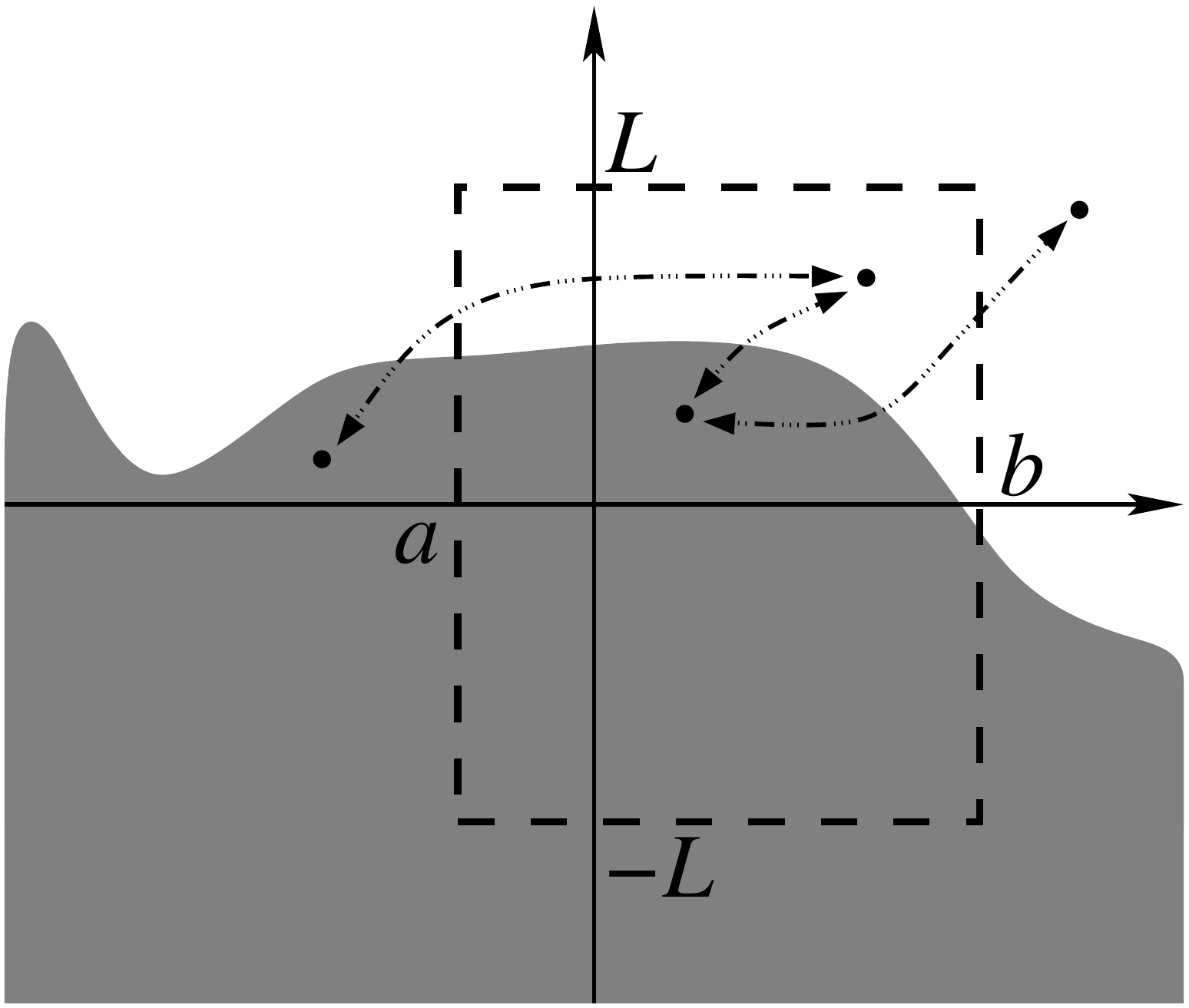}
	\caption{\sl The interactions contributing to the fractional perimeter.}
			\label{DIFI1oinfkerfmD}
\end{figure}

The mathematical notation used in this paper goes as follows. Given~$s\in(0,1)$, $a<b$ 
and~$u_0\in C(\R)$,
we say that~$u\in L^\infty_{\rm loc}(\R)$ is an $s$-minimal graph in~$(a,b)$ if:
\begin{itemize}
\item $u(x)=u_0(x)$ for every~$x\in\R\setminus(a,b)$,
\item for every~$L>\|u\|_{L^\infty((a,b))}$ and every measurable function~$v:\R\to\R$ such that~$|v(x)|\le L$
for all~$x\in(a,b)$ we have that
$$ \Per_s(E_u,\Omega_L)\leq\Per_s(E_v,\Omega_L)
$$
where
\begin{eqnarray*}
&& E_u:=\big\{ (x,y)\in\R^2 {\mbox{ s.t. }} y<u(x)\big\},
\\&&E_v:=\big\{ (x,y)\in\R^2 {\mbox{ s.t. }} y<v(x)\big\},\\
&&\Omega_L:=(a,b)\times(-L,L),\\
&& \Per_s(E,\Omega):={\mathcal{L}}_s(E\cap\Omega,E^c\cap\Omega)+
{\mathcal{L}}_s(E\cap\Omega,E^c\cap\Omega^c)+
{\mathcal{L}}_s(E\cap\Omega^c,E^c\cap\Omega),\\
&& \Omega^c:=\R^2\setminus\Omega\\
{\mbox{and }}&&{\mathcal{L}}_s(A,B):=\iint_{A\times B}\frac{dX \,dY}{|X-Y|^{2+s}}.
\end{eqnarray*}
\end{itemize}
See Figure~\ref{DIFI1oinfkerfmD} for a sketch of the interactions detected by the definition above.
See also~\cite{MR4279395} for further information about nonlocal minimal graphs
and in particular Theorem~1.11 there
for existence and uniqueness results for nonlocal minimal graphs.

We also recall that, while the minimizers of the classical perimeter have vanishing mean curvature,
a minimizer~$E$ in~$\Omega$ of the fractional perimeter has vanishing fractional mean curvature, in the sense that,
for every~$P\in\Omega\cap(\partial E)$,
$$ \int_{\R^2}\frac{\chi_{\R^2\setminus E}(X)-\chi_{E}(X)}{|X-P|^{2+s}}\,dX=0.$$
Here above and in the sequel, the singular integral is intended in the Cauchy principal value
sense and the above equation reads in the sense of viscosity, see~\cite[Theorem~5.1]{MR2675483} for full details.\medskip

Numerical examples of the stickiness properties have recently been
provided by~\cite{MR3982031, MR4294645}.
The theory of nonlocal minimal surfaces also presents a number of interesting offsprings, such as regularity
theory~\cite{MR3090533, MR3107529, MR3331523, MR3680376, MR3798717, MR3934589, MR3981295, MR4050198, MR4116635}, isoperimetric problems and the study of constant mean curvature surfaces~\cite{MR2799577,MR3322379, MR3412379, MR3485130, MR3744919, MR3836150, MR3881478},
geometric evolution problems~\cite{MR2487027, MR3156889, MR3713894, MR3778164, MR3951024}, front propagation problems~\cite{MR2564467}, phase transition problems~\cite{MR3133422}, etc.

\medskip

Having formalized the setting in which we work, we can now state our result
about the stickiness phenomenon in the antisymmetric setting as an arbitrarily small perturbation of the flat case:

\begin{theorem}\label{MAIN:T}
Let~$\epsilon_0>0$, $\bar{C}>0$, $d>d_0>0$ and~$h>0$, with
\begin{equation}\label{KS-GEOP} \frac{h^{2+s}}{(h^2+2)^{\frac{2+s}{2}}}>
\frac{2\bar{C} }{(1+s){\left(d-d_0\right)^{1+s}}}
.\end{equation}
Then, there exist~$C>0$ and~$\eta_0\in(0,1)$, depending only
on~$s$, $\epsilon_0$, $\bar{C}$, $d$, $d_0$ and~$h$, such that if~$\eta\in(0,\eta_0]$
the following claim holds true.

Let~$u_0\in C(\R)$.
Assume that
\begin{equation}\label{LA8765S-OFDFetahr}
u_0(x)=0 \quad{\mbox{ for all }}x\in(d,d+h),\end{equation}
\begin{equation}\label{SM12425G-2}
u_0(x)\le0\quad{\mbox{ for all }}x\in(d,+\infty),
\end{equation}
\begin{equation}\label{LA8765S-OFDFeta}
\|u_0\|_{L^\infty(\R)}\le \bar{C}\eta,\end{equation}
\begin{equation}\label{LA8765S-OFDF}
u_0(x)=-u_0(-x)\quad{\mbox{ for all }}x\in(d,+\infty)
\end{equation}
and
\begin{equation}\label{LA8765S-OFDF-CSTA}
\int_{-\infty}^{-d-h}\frac{u_0(x)}{|d_0-d-x|^{2+s}}\,dx\ge \eta.\end{equation}
Let~$u$ be the $s$-minimal graph with~$u=u_0$ in~$\R\setminus(-d,d)$.

Then, 
\begin{equation}\label{SMG-3-BISla}
u(x)=-u(-x)\quad{\mbox{ for all }}x\in(-\infty,+\infty),
\end{equation}
\begin{equation}\label{SMG-4-BISla}
u(x)\ge0\quad{\mbox{ for all }}x\in(-\infty,0],
\end{equation}
\begin{equation}\label{SMG-5-BISla}
u(x)\le0\quad{\mbox{ for all }}x\in[0,+\infty)
\end{equation}
and, for every~$x\in(-d,-d+d_0)$,
\begin{equation}\label{SMG-BISlaX} u(x)\ge \frac{\eta^{\frac{2+\epsilon_0}{1-s}}}{C}.\end{equation}
\end{theorem}

\begin{figure}
		\centering
		\includegraphics[width=.9\linewidth]{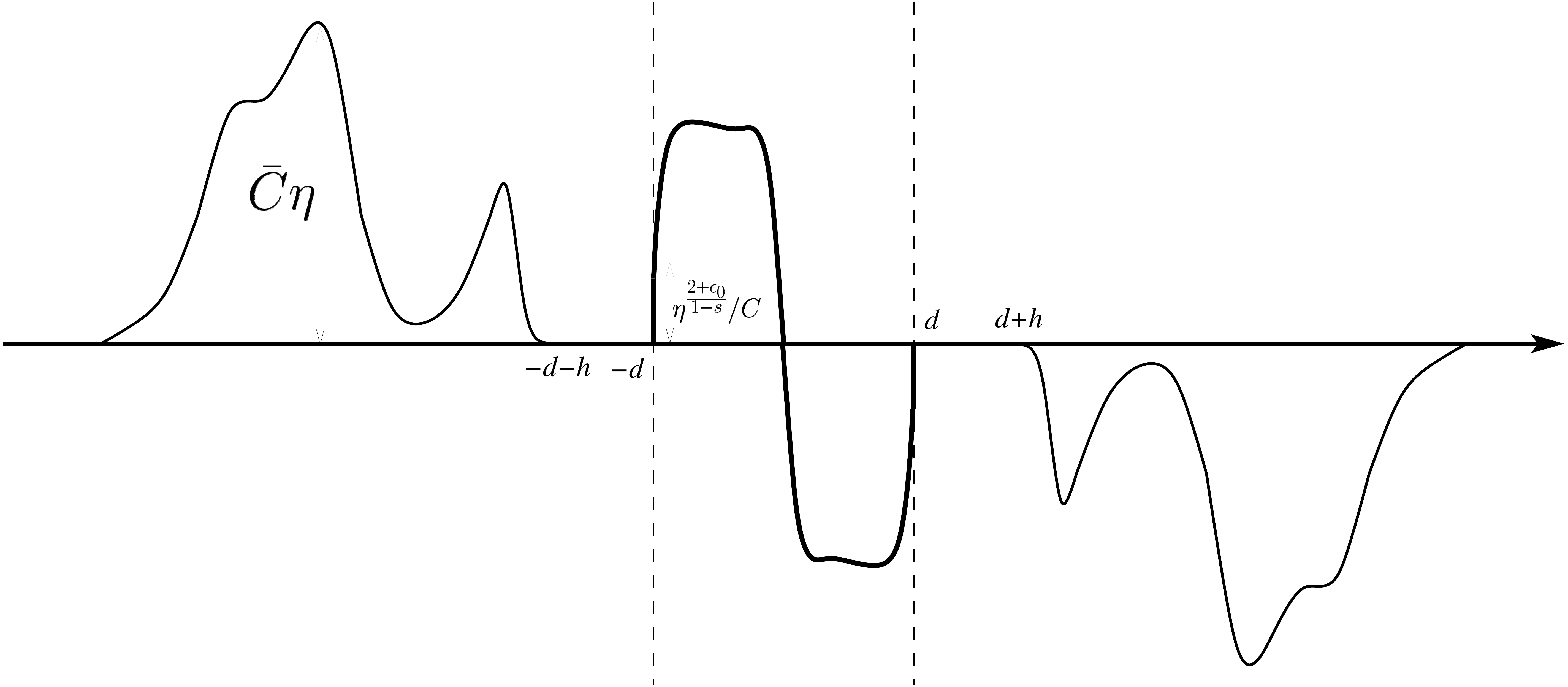}
	\caption{\sl The antisymmetric stickiness phenomenon, as given in Theorem~\ref{MAIN:T}.}
	\label{3wqf4DIFI2}
\end{figure}

See Figure~\ref{3wqf4DIFI2} for a sketch of the geometric scenario described in Theorem~\ref{MAIN:T}.
Interestingly, the antisymmetric stickiness power~$\frac{2+\epsilon_0}{1-s}$ in~\eqref{SMG-BISlaX}
is the same as the one detected in~\cite[Theorem~1.4]{MR3596708} for the even symmetric case
(hence, the antisymmetric case appears to be also quantitatively in agreement with the
even symmetric configurations, up to multiplicative constants).

Notice also that condition~\eqref{KS-GEOP} simply states that~$d$ is sufficiently large with respect to other structural parameters.
For instance, one can take
\begin{eqnarray*}
&&\bar{C}:=2^{2+s}(1+s),\qquad d:=3^{\frac{2+s}{2(1+s)}}2^{\frac{3+s}{1+s}}
+2
\qquad d_0:=1 \quad{\mbox{ and }} \quad h:=1.
\end{eqnarray*}
Notice that with this choice of parameters, condition~\eqref{KS-GEOP}
is satisfied.
We also define
$$ d_1:=2\left(\left(\frac{10}{9}\right)^{\frac1{1+s}}-1\right)
\quad{\mbox{ and }}\quad
d_2:=2\left(10^{\frac1{1+s}}-1\right),$$
and we take~$u_0\in C(\R)$ 
satisfying~\eqref{LA8765S-OFDFetahr}, \eqref{SM12425G-2},
\eqref{LA8765S-OFDFeta}
and~\eqref{LA8765S-OFDF}, and such
that
$$ u_0\ge\bar{C}\eta\chi_{(-d-h-d_2,-d-h-d_1)}\quad {\mbox{ in }}(-\infty,-d).$$
With this choice of~$u_0$ we have that~\eqref{LA8765S-OFDF-CSTA} is also
satisfied. Indeed, using the change of variable~$z:=d_0-d-x$,
\begin{eqnarray*}&&
\int_{-\infty}^{-d-h}\frac{u_0(x)}{|d_0-d-x|^{2+s}}\,dx
\ge \int_{-d-h-d_2}^{-d-h-d_1}\frac{u_0(x)}{|d_0-d-x|^{2+s}}\,dx\\&&\qquad\qquad
\ge \int_{-d-h-d_2}^{-d-h-d_1}\frac{\bar{C}\eta}{|d_0-d-x|^{2+s}}\,dx
=\bar{C}\eta \int_{d_0+h+d_1}^{d_0+h+d_2}\frac{dz}{z^{2+s}}\\&&\qquad\qquad=
\frac{\bar{C}\eta}{1+s}\left[\frac1{(d_0+h+d_1)^{1+s}} -\frac1{(d_0+h+d_2)^{1+s}}
\right]\\
&&\qquad\qquad=
2^{2+s}\eta \left(\frac1{2^{1+s}\frac{10}{9}} -\frac1{2^{1+s} 10}
\right)=2 \eta\frac{8}{10} >\eta,
\end{eqnarray*}
as desired.

We also point out that condition~\eqref{LA8765S-OFDFetahr} can be relaxed
and is taken here mostly to emphasize that ``instability''
(as embodied by the stickiness phenomenon)
arises in the case of nonlocal minimal surfaces even when the external datum is,
near the boundary of the domain,
as regular, as ``stable'' and ``as flat'' as one wishes: more specifically, the stickiness that
we detect is not due to possible
oscillatory behaviors of the external datum near the domain of reference, but is merely the 
outcome of long-distance interactions.

As for condition~\eqref{LA8765S-OFDF-CSTA}, its meaning is, roughly speaking, that the 
external datum presents a
``positive bump'' at the left (which corresponds to a ``negative bump'' at the right).
The role of the parameter~$\eta$ is to give an explicit quantification of this bump: this 
quantification
is given in an integral form, since we want to weigh the ``size'' of the bump (given by the 
numerator
of the integrand on the left hand side of~\eqref{LA8765S-OFDF-CSTA}) against the 
interaction kernel,
because this is, in a sense, the ``long-distance force'' exerted by the bump on the nonlocal 
minimal graph.
\medskip

The gist of the proof of Theorem~\ref{MAIN:T} is to employ the auxiliary function constructed 
in~\cite[Corollary~7.2]{MR3596708}.
This barrier will be suitably scaled and placed near the left end of the 
reference domain
to show that the nonlocal minimal graph under consideration must be ``lifted up''.
{F}rom a quantitative point of view, this localized barrier
presents a nonlocal mean curvature which
has possibly a ``wrong sign'' somewhere, but this sign discrepancy is controlled
by a small quantity (roughly speaking, playing the role of~$\eta$ in the statement
of Theorem~\ref{MAIN:T}, up to a convenient choice
of a multiplicative constant).
With this respect, the external bump in condition~\eqref{LA8765S-OFDF-CSTA}
compensates the ``small sign discrepancy''
of this localized barrier.
That is, roughly speaking, the positive bump on the left will produce an
{\em advantageous term}, while the negative bump on the right and the modification
needed to lower the barrier in the negativity region of~$u$
will produce a {\em disadvantageous
term of the same order}. Thus, playing around with constants, one detects a natural
structure assuring that the advantageous term is greater than the sum of the initial
error on the nonlocal mean curvature, the contribution of the negative bump and the
error coming from the barrier modification.

In this procedure, however, an additional difficulty arises, since in our framework the
solution will cross the horizontal axis. For this reason, the previous barrier (which is
designed to be a small upward bump) cannot be exploited as it is and needs to be
modified to stay below the ``expected'' negative regions of the minimizer.

To this end, we first need to identify a ``safe region near the boundary'' for the minimizer,
that is an interval of well-determined length in which one can be sure that the solution
is positive: once this is accomplished, it will be possible
to slide from below the appropriate barrier whose positive portion occurs precisely in
the above safe region.
A graphical sketch of such a barrier will be given in Figure~\ref{21421533wqf4DIFI2}.
\medskip

In view of these comments, to identify the safe region near
the boundary and perform the proof of
Theorem~\ref{MAIN:T}, we establish a result of general interest, which can be seen
as a ``maximum principle for
antisymmetric $s$-minimal graphs''. 
For our purposes, the safe region
will correspond to the positive values of~$u$, which, thanks to this antisymmetric
maximum principle, extends to the whole interval~$(-d,0)$.

Roughly speaking, the classical maximum principle for
$s$-minimal graphs (see e.g.~\cite{MR2675483, MR3516886}) states that if the external
data of
an $s$-minimal graph are positive (or negative) then so is the $s$-minimal graph.
Here we deal with an antisymmetric configuration, hence it is not possible for the 
external data of an $s$-minimal graph to have a sign (except in the trivial case of
identically zero datum,
which
produces the identically zero $s$-minimizer). Hence, the natural counterpart of the
maximum principle in
the antisymmetric framework is to assume that the external datum ``on one side'' has
a sign, thus forcing
the external datum on the other side to have the opposite sign.
Under this assumption, 
we show that the corresponding $s$-minimal graph maintains the antisymmetry
and sign properties of the datum, according
to the following result:

\begin{theorem}\label{SMG-T}
Let~$d$, $h>0$. Let~$u:\R\to\R$ be an $s$-minimal graph in~$(-d,d)$, with~$u\in C(-\infty,-d)\cap C^{1,\frac{1+s}2}(-d-h,-d)$.
Assume that
\begin{equation}\label{SMG-1TH}
u(x)=-u(-x)\quad{\mbox{ for all }}x\in(d,+\infty)
\end{equation}
and that
\begin{equation}\label{SMG-2TH}
u(x)\le0\quad{\mbox{ for all }}x\in(d,+\infty).
\end{equation}

Then,
\begin{equation}\label{SMG-3}
u(x)=-u(-x)\quad{\mbox{ for all }}x\in(-\infty,+\infty),
\end{equation}
\begin{equation}\label{920irpkfeitk-320jtjoyj-01}
u(x)\ge0\quad{\mbox{ for all }}x\in(-\infty,0]
\end{equation}
and
\begin{equation}\label{920irpkfeitk-320jtjoyj-02}
u(x)\le0\quad{\mbox{ for all }}x\in[0,+\infty).
\end{equation}
\end{theorem}

\begin{figure}
		\centering
		\includegraphics[width=.9\linewidth]{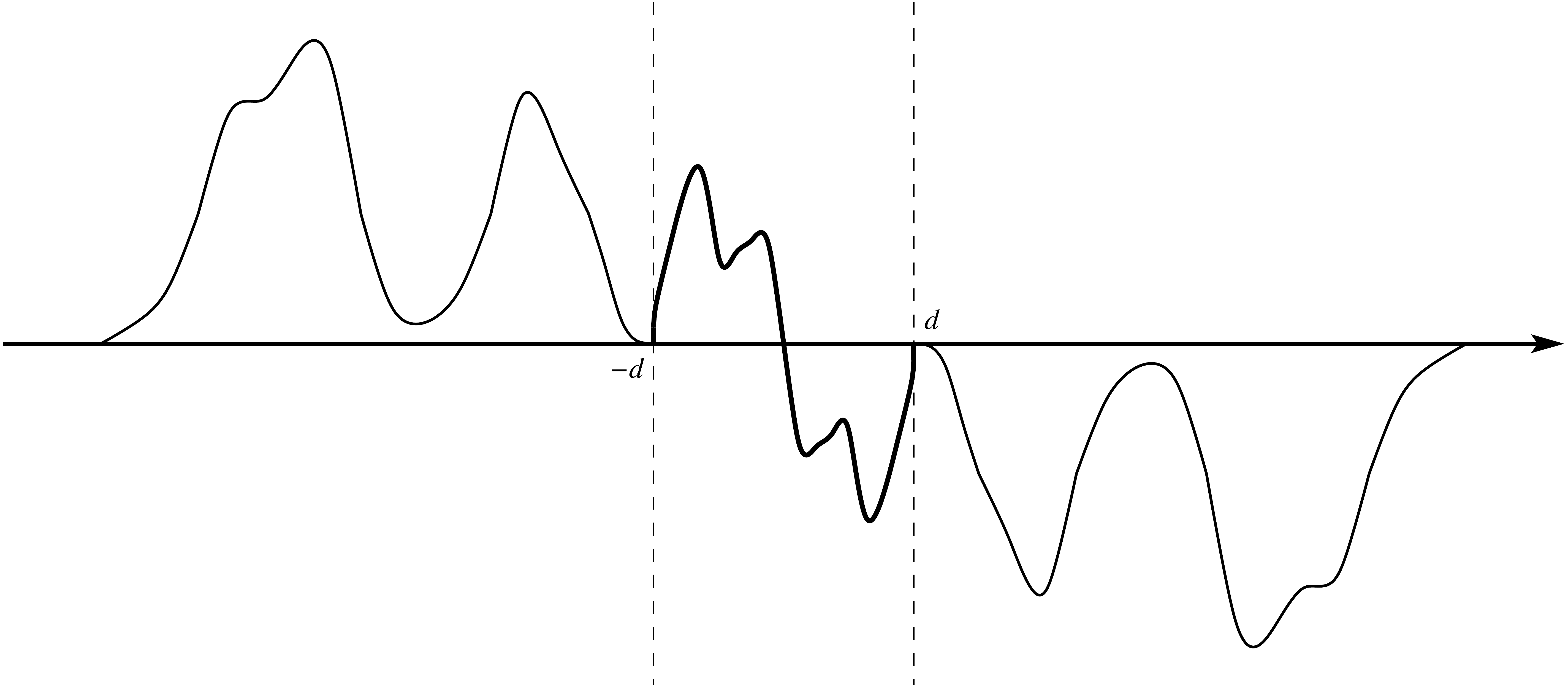}
	\caption{\sl The maximum principle for antisymmetric nonlocal minimal graphs, as given in Theorem~\ref{SMG-T}.}
	\label{3DIFI2}
\end{figure}

See Figure~\ref{3DIFI2} for a sketch of the configuration detected in Theorem~\ref{SMG-T}.

The proofs of Theorems~\ref{MAIN:T} and~\ref{SMG-T} are contained in Sections~\ref{SEC-PF1}
and~\ref{SEC-PF2}, respectively. As a matter of fact,
Theorem~\ref{SMG-T} will follow from a more general statement valid for supersolutions, given in the forthcoming Lemma~\ref{SUSO}.

Section~\ref{SNEW} contains
a discussion about the proof of the antisymmetric maximum principle,
compared with the existing theory focused on the antisymmetric maximum principle for the fractional Laplacian
and linear nonlocal operators.

\section{Proof of Theorem~\ref{MAIN:T}}\label{SEC-PF1}

We give here the proof of Theorem~\ref{MAIN:T}.
At this stage, we freely use the antisymmetric maximum principle in
Theorem~\ref{SMG-T}, whose proof is postponed to Section~\ref{SEC-PF2}.

We also use the notation~$X=(x,y)$ to denote points in~$\R^2$, with~$x$, $y\in\R$.

\begin{figure}
		\centering
		\includegraphics[width=.9\linewidth]{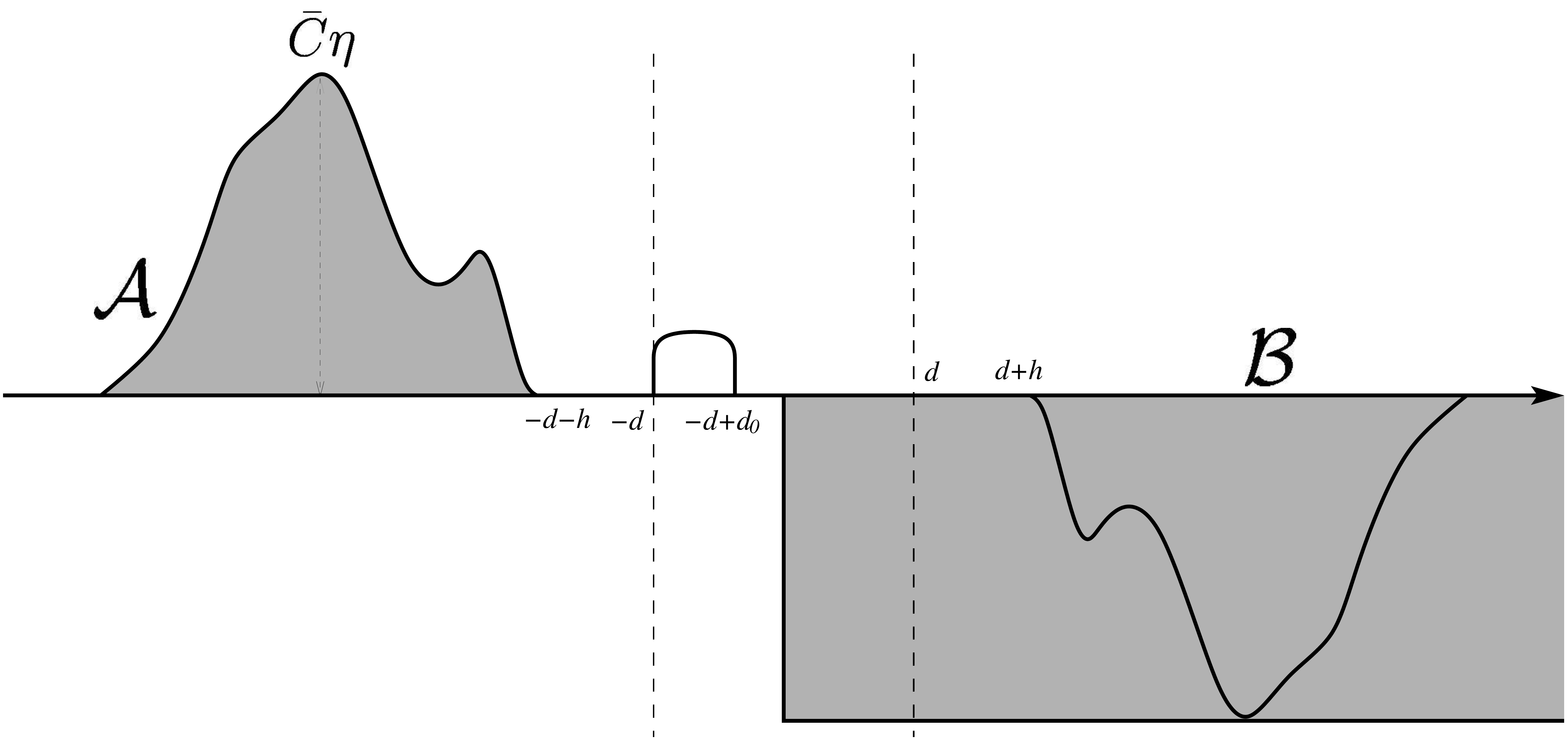}
	\caption{\sl The sets~${\mathcal{A}}$ and~${\mathcal{B}}$ in the proof of Theorem~\ref{MAIN:T}.}
			\label{DIFI9qqywieohgk-10}
\end{figure}

The proof of Theorem~\ref{MAIN:T} consists of five steps.\medskip

\noindent {\bf{Step 1:}} {\em Construction of the barrier~$U_\delta$.}
Let~$\epsilon_0>0$.
Given~$\delta\in(0,1)$ sufficiently small, we apply~\cite[Corollary~7.2]{MR3596708}
and we find a set~$H_\delta\subset\R^2$
that contains the halfplane~$\R\times(-\infty,0)$
and is contained in the halfplane~$\R\times(-\infty,\delta)$, such that
\begin{eqnarray*}&&
\left\{(x,y)\in H_\delta {\mbox{ s.t. }} x\in(-\infty,-d)\cup\left(-d+d_0,+\infty\right)
\right\}\\&&\qquad=  \left\{(x,y)\in \R^2 {\mbox{ s.t. }} 
x\in(-\infty,-d)\cup\left(-d+d_0,+\infty\right) {\mbox{ and }}
y\in(-\infty,0)\right\},\end{eqnarray*}
with
\begin{equation}\label{SNDomwrfeRGHSo}
\begin{split}& \left\{(x,y)\in H_\delta {\mbox{ s.t. }} x\in\left(-d,-d+d_0\right)\right\}\\&\qquad\supseteq  \left\{(x,y)\in \R^2 {\mbox{ s.t. }} x\in
\left(-d,-d+d_0\right)\times\left(-\infty,\delta^{\frac{2+\epsilon_0}{1-s}}\right)\right\}\end{split}\end{equation}
and satisfying, for each~$P=(p,q)\in \partial H_\delta$ with~$p\in\left(-d,-d+d_0\right)$,
\begin{equation}\label{USJN-iaA} \int_{\R^2}\frac{\chi_{\R^2\setminus H_\delta}(X)-\chi_{H_\delta}(X)}{|X-P|^{2+s}}\,dX\le C\delta,\end{equation}
for some constant~$C>0$ depending only on~$s$, $\epsilon_0$ and~$d$.

\begin{figure}
		\centering
		\includegraphics[width=.9\linewidth]{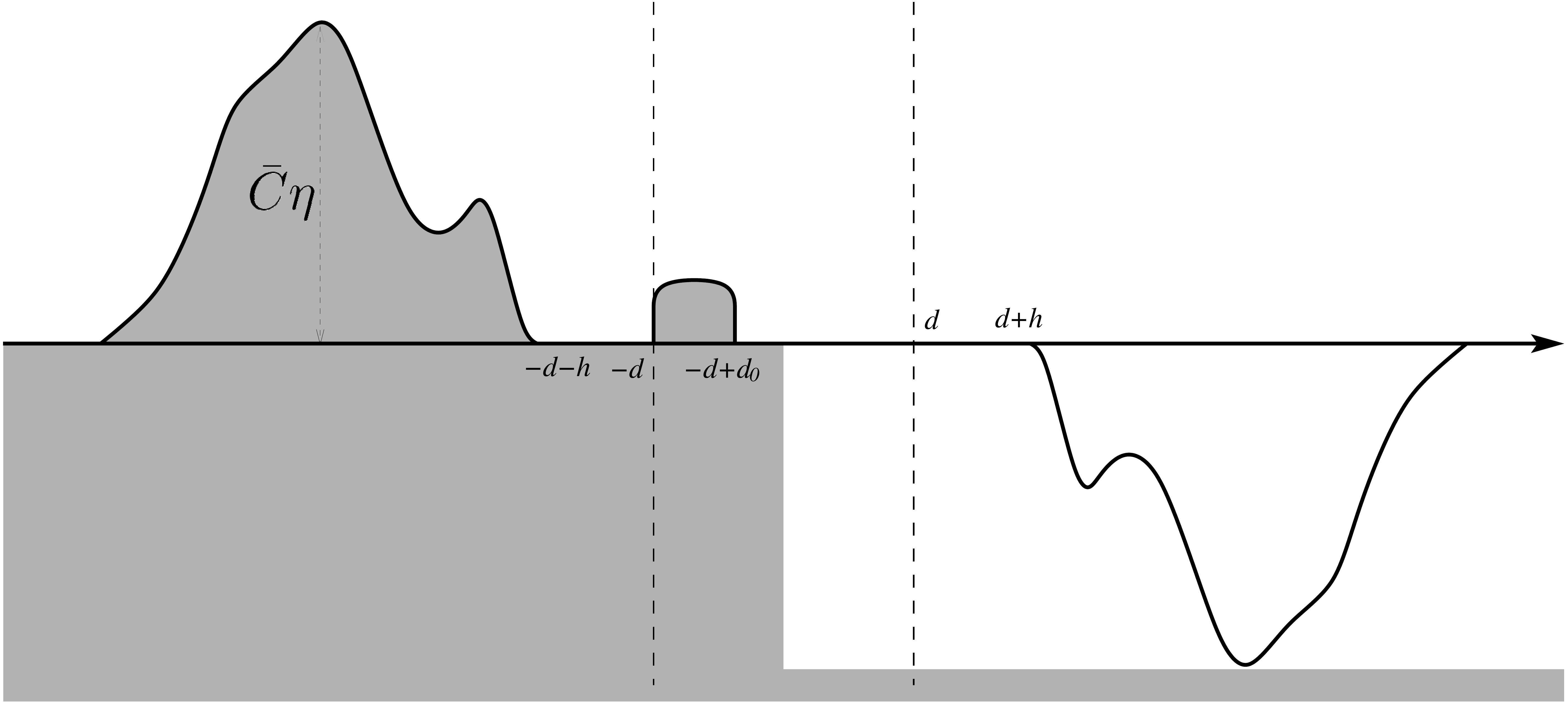}
	\caption{\sl The set~$U_\delta$ used in the proof of Theorem~\ref{MAIN:T}.}
	\label{21421533wqf4DIFI2}
\end{figure}

The strategy will be to take
\begin{equation}\label{MAIN:T:ETADEL}
\eta:=C_\star\delta\end{equation}
in the statement of Theorem~\ref{MAIN:T}, with~$C_\star>0$ to be sufficiently large.
We let
\begin{eqnarray*}&& {\mathcal{A}}:=\big\{(x,y)\in \R^2 {\mbox{ s.t. }} x\in(-\infty,-d) {\mbox{ and }} 0<y<u_0(x)\big\},
\\&&{\mathcal{B}}:=\big\{(x,y)\in \R^2 {\mbox{ s.t. }} x\in(0,+\infty)
 {\mbox{ and }} -\bar{C}\eta<y<0\big\}\\
{\mbox{and}}\quad&& U_\delta:=(H_\delta\cup {\mathcal{A}})
\setminus {\mathcal{B}},\end{eqnarray*}
see Figure~\ref{DIFI9qqywieohgk-10} to visualize the sets~${\mathcal{A}}$ and~${\mathcal{B}}$
and Figure~\ref{21421533wqf4DIFI2}
to visualize the set~$U_\delta$.
\medskip

\noindent {\bf{Step 2:}} {\em Estimates for some integrals appearing in the nonlocal mean curvature
of~$U_\delta$.}
We observe that~${\mathcal{A}}\cap H_\delta= \varnothing$ and~${\mathcal{B}}
\subset H_\delta$. Consequently, 
\begin{eqnarray*}
&&\chi_{U_\delta}=\chi_{H_\delta}+\chi_{{\mathcal{A}}}-\chi_{{\mathcal{B}}}\\
{\mbox{and }}&&
\chi_{\R^2\setminus U_\delta}=\chi_{\R^2\setminus(H_\delta\cup{\mathcal{A}})
}+\chi_{{\mathcal{B}}}\le \chi_{\R^2\setminus H_\delta 
}+\chi_{{\mathcal{B}}}.
\end{eqnarray*}
{F}rom this and the inequality
for~$H_\delta$ in~\eqref{USJN-iaA}, we have that,
for all~$P=(p,q)\in \partial H_\delta$ with~$p\in\left(-d,-d+d_0\right)$,
\begin{equation}\label{M24R5YN-OAM78}\begin{split}&
\int_{\R^2}\frac{\chi_{\R^2\setminus 
U_\delta}(X)-\chi_{U_\delta}(X)}{|X-P|^{2+s}}\,dX\\&\qquad\le
\int_{\R^2}\frac{\chi_{\R^2\setminus H_\delta 
}(X)+\chi_{{\mathcal{B}}}(X)-\big(
\chi_{H_\delta}(X)+\chi_{{\mathcal{A}}}(X)-\chi_{{\mathcal{B}}(X)}\big) }{|X-P|^{2+s}}\,dX
\\&\qquad\le C\delta
+2\int_{{\mathcal{B}}}\frac{dX}{|X-P|^{2+s}}-\int_{{\mathcal{A}}}\frac{dX}{|X-P|^{2+s}}.
\end{split}\end{equation}

We note that, for all~$P=(p,q)\in \partial H_\delta$ with~$p\in\left(-d,-d+d_0\right)$,
\begin{eqnarray*}
&&\int_{{\mathcal{B}}}\frac{dX}{|X-P|^{2+s}}\le\int_{{\mathcal{B}}}\frac{dX}{|x-p|^{2+s}}=
\bar{C}\eta\int_{0}^{+\infty}\frac{dx}{|x-p|^{2+s}}=\bar{C}\eta\int_{0}^{+\infty}\frac{dx}{(x-p)^{2+s}}\\&&\qquad\le
\bar{C}\eta\int_{0}^{+\infty}\frac{dx}{\left(x+d-d_0\right)^{2+s}}
=\frac{\bar{C}\eta}{(1+s){\left(d-d_0\right)^{1+s}}}
=\frac{\bar{C} C_\star\delta}{(1+s){\left(d-d_0\right)^{1+s}}}.
\end{eqnarray*}
This and~\eqref{M24R5YN-OAM78} yield that
\begin{equation}\label{M24R5YN-OAM78-0osSKMx0e}
\int_{\R^2}\frac{\chi_{\R^2\setminus U_\delta}(X)-\chi_{U_\delta}(X)}{
|X-P|^{2+s}}\,dX\le C\delta
+\frac{2\bar{C} C_\star\delta}{(1+s){\left(d-d_0\right)^{1+s}}}
-\int_{{\mathcal{A}}}\frac{dX}{|X-P|^{2+s}}.
\end{equation}\medskip

\noindent{\bf{Step 3:}} {\em Estimate for the third term in~\eqref{M24R5YN-OAM78-0osSKMx0e}.}
Let us now estimate the last integral in~\eqref{M24R5YN-OAM78-0osSKMx0e}.
By~\eqref{LA8765S-OFDFeta}, if~$X=(x,y)\in{{\mathcal{A}}}$ and~$P=(p,q)\in \partial H_\delta$ with~$p\in\left(-d,-d+d_0\right)$,
$$ |y- q|\le|y|+|q|\le1+\delta\le 2.$$
Moreover, owing to the assumptions on~$u$ in~\eqref{LA8765S-OFDFetahr} and~\eqref{LA8765S-OFDF}, for all~$x\in(-d-h,-d)$ we have that~$
u_0(x)=-u_0(-x)=0$. 

Accordingly, if~$X=(x,y)\in{{\mathcal{A}}}$ then~$x\le-d-h$. As a result,
if~$X=(x,y)\in{{\mathcal{A}}}$ and~$P=(p,q)\in \partial H_\delta$ with~$p\in\left(-d,-d+d_0\right)$,
\begin{equation*}\begin{split}& |x- p|\ge\min\big\{ |x-p|,\,|x+p|\big\}
\ge\min\big\{ p-x,\,-p-x\big\}\\&\qquad\ge\min\left\{ -d-(-d-h),\, 
d-d_0-(-d-h)\right\}=
\min\left\{h, 2d-d_0+h\right\}=h.\end{split}\end{equation*}

{F}rom these remarks we obtain that
\begin{eqnarray*}&&|X-P|=|x-p|\sqrt{1+\frac{|y-q|^2}{|x-p|^2}}
\le|x-p|\sqrt{1+\frac{2}{h^2}}=|x-p|\sqrt{\frac{h^2+2}{h^2}}.
\end{eqnarray*}
Hence,
\begin{equation}\label{NBSinSUJHNDgYDwqr4e} \begin{split}&\int_{{\mathcal{A}}}\frac{dX}{|X-P|^{2+s}}\ge\frac{h^{2+s}}{(h^2+2)^{\frac{2+s}{2}}}\int_{{\mathcal{A}}}\frac{dX}{|x-p|^{2+s}}\\&\qquad\qquad=\frac{h^{2+s}}{(h^2+2)^{\frac{2+s}{2}}}
\int_{-\infty}^{-d}\frac{u_0(x)}{|x-p|^{2+s}}\,dx=\frac{h^{2+s}}{(h^2+2)^{\frac{2+s}{2}}}\int_{-\infty}^{-d-h}\frac{u_0(x)}{|x-p|^{2+s}}\,dx.\end{split}\end{equation}

Since, if~$p\in\left(-d,-d+d_0\right)$ and~$x\in(-\infty,-d-h)$,
$$ |x-p|=p-x\le-d+d_0-x=|d_0-d-x|,
$$
we deduce from~\eqref{NBSinSUJHNDgYDwqr4e} that
$$ \int_{{\mathcal{A}}}\frac{dX}{|X-P|^{2+s}}\ge\frac{h^{2+s}}{(h^2+2)^{\frac{2+s}{2}}}
\int_{-\infty}^{-d-h}\frac{u_0(x)}{|d_0-d-x|^{2+s}}\,dx.$$
This, the integral assumption
on~$u_0$ in~\eqref{LA8765S-OFDF-CSTA} and the definition of~$\eta$ in~\eqref{MAIN:T:ETADEL} give that
$$ \int_{{\mathcal{A}}}\frac{dX}{|X-P|^{2+s}}\ge \frac{h^{2+s}\eta}{(h^2+2)^{\frac{2+s}{2}}}=\frac{C_\star h^{2+s}\delta}{(h^2+2)^{\frac{2+s}{2}}}.$$
Plugging this information into~\eqref{M24R5YN-OAM78-0osSKMx0e} we infer that
\begin{equation}\label{KS-GEOPXM}
\int_{\R^2}\frac{\chi_{\R^2\setminus U_\delta}(X)-\chi_{U_\delta}(X)}{|X-P|^{2+s}}\,dX\le C\delta
+\frac{2\bar{C} C_\star\delta}{(1+s){\left(d-d_0\right)^{1+s}}}-
\frac{C_\star h^{2+s}\delta}{(h^2+2)^{\frac{2+s}{2}}}.
\end{equation}
\medskip

\noindent{\bf{Step 4:}} {\em Conclusion of the proof that~$U_\delta$ is a subsolution.}
Now we define
$$ \vartheta:=
\frac{h^{2+s}}{(h^2+2)^{\frac{2+s}{2}}}-
\frac{2\bar{C} }{(1+s){\left(d-d_0\right)^{1+s}}}
$$
and we observe that~$\vartheta>0$, owing to the relation
between parameters in~\eqref{KS-GEOP}.
Then, we deduce from~\eqref{KS-GEOPXM} that
\begin{equation}\label{KS-GEOPXM2}
\int_{\R^2}\frac{\chi_{\R^2\setminus U_\delta}(X)-\chi_{U_\delta}(X)}{|X-P|^{2+s}}\,dX\le C\delta
-C_\star\vartheta\delta\le-\frac{C_\star\vartheta\delta}2<0,
\end{equation}
as long as~$C_\star$ is large enough.\medskip

\noindent{\bf{Step 5:}} {\em Sliding method and conclusion of the proof of Theorem~\ref{MAIN:T}.}
We can now use~$U_\delta$ as a barrier for the sliding method. Specifically, by~\eqref{LA8765S-OFDFeta} and the maximum principle in~\cite{MR2675483}, we know that~$|u(x)|\le\bar{C}\eta<1$ for all~$x\in\R$ and therefore a downwards translation
of~$U_\delta$ with magnitude~$1$, that we denote by~$U_\delta-(0,1)$,
is completely contained in~$E_u:=\big\{ (x,y)\in\R^2 {\mbox{ s.t. }} y<u(x)\big\}$.

We then slide  up till we reach a touching point: given~$\tau\in(0,1]$, we notice that by construction no touching can occur between~$E_u$ and~$U_\delta-(0,\tau)$ at points with abscissa in~$\R\setminus[-d,0]$.
But these touching points cannot occur in~$\left(-d+d_0,0\right]$ either, thanks to the sign assumption in~\eqref{SMG-4}.
And they cannot occur in~$\left[-d,-d+d_0\right]$, thanks to~\eqref{KS-GEOPXM2}
and the maximum principle in~\cite{MR2675483} (see also~\cite[Theorem~1.4]{MR4104542}).

As a consequence, we have that~$U_\delta\subseteq E_u$.
By virtue of~\eqref{SNDomwrfeRGHSo} and~\eqref{MAIN:T:ETADEL}, this establishes~\eqref{SMG-BISlaX}.

Also, the claims in~\eqref{SMG-3-BISla}, \eqref{SMG-4-BISla} and~\eqref{SMG-5-BISla}
follow from~\eqref{SMG-3}, \eqref{SMG-4} and~\eqref{SMG-5}.
\hfill$\Box$

\section{Proof of Theorem~\ref{SMG-T}}\label{SEC-PF2}

We now provide a more general antisymmetric maximum principle
valid for supersolutions, from which Theorem~\ref{SMG-T} will plainly follow by using also other results already available in the literature.

\begin{lemma}\label{SUSO}
Let~$d>0$. Let~$u:\R\to\R$, with~$u\in C([-d,0])\cap C^{1,\alpha}((-d,0])$
for some~$\alpha\in\left(\frac{1+s}2,1\right)$.

Let also
$$E:=\big\{(x,y)\in\R^2 {\mbox{ s.t. }}y<u(x)\big\}.$$
Assume that, for all~$P=(p,u(p))$ with~$p\in(-d,0]$, we have
\begin{equation}\label{S:10} \int_{\R^2}\frac{\chi_{\R^2\setminus E}(X)-\chi_{E}(X)}{|X-P|^{2+s}}\,dX\ge0.\end{equation}
Suppose also that
\begin{equation}\label{SECCA}
u(-d)>\min_{[-d,0]}u,
\end{equation}
that
\begin{equation}\label{SMG-1}
u(x)=-u(-x)\quad{\mbox{ for all }}x\in(0,+\infty)
\end{equation}
and that
\begin{equation}\label{SMG-2}
u(x)\le0\quad{\mbox{ for all }}x\in(d,+\infty).
\end{equation}

Then,
\begin{equation}\label{SMG-4}
u(x)\ge0\quad{\mbox{ for all }}x\in(-\infty,0]
\end{equation}
and
\begin{equation}\label{SMG-5}
u(x)\le0\quad{\mbox{ for all }}x\in[0,+\infty).
\end{equation}
\end{lemma}

\begin{proof} It suffices to show that~$u\ge0$ in~$(-\infty,0]$, that is~\eqref{SMG-4},
since~\eqref{SMG-5} would then follow from
the odd symmetry of~$u$ given by~\eqref{SMG-1}
and the nonnegativity of~$u$ in~$(-\infty,0]$ in~\eqref{SMG-4}.

To prove~\eqref{SMG-4} we argue by contradiction and assume that~\eqref{SMG-4}
is violated. Note that
\begin{equation}\label{SMG-8}
{\mbox{if~$x\in(-\infty,-d)$ then~$u(x)=-u(-x)\ge0$,}}\end{equation}
thanks to the assumptions on~$u$ in~\eqref{SMG-1} and~\eqref{SMG-2}.

This, the continuity of~$u$ in~$[-d,0]$
and our contradictory assumption give that there exists~$p\in [-d,0]$ such that
\begin{equation}\label{UYP} u(p)=\min_{[-d,0]} u<0.\end{equation}

We claim that
\begin{equation}\label{SMG-L}
p\in(-d,0).\end{equation}
Indeed,
since~$u(0)=-u(0)$ by the
odd symmetric of~$u$ in~\eqref{SMG-1}, we have that~$u(0)=0$ and accordingly~$p\ne0$. 
This and~\eqref{SECCA} lead to~\eqref{SMG-L}, as desired.

By~\eqref{SMG-L}, we know that
the supersolution condition~\eqref{S:10} holds true for~$P:=(p,u(p))$.

We will now reach a contradiction by showing that~\eqref{S:10} is violated, since ``the set~$E$ is too large''.
To this end, we consider the sets
\begin{equation}\label{ieuryurfbsfbsjffjfghfweeu09998ye73}\begin{split}
& F_-:=\big\{x\in\R {\mbox{ s.t. }} u(x)<u(p)\big\}\\
{\mbox{and }}\quad &F_+:=\big\{x\in\R {\mbox{ s.t. }} u(x)> -u(p)\big\},
\end{split}\end{equation}
see Figure~\ref{DIFI2} (notice in particular that
the sets~$F_-$ and~$F_+$ are given by the bold intervals along the
horizontal axis in  Figure~\ref{DIFI2}).

\begin{figure}
		\centering
		\includegraphics[width=.9\linewidth]{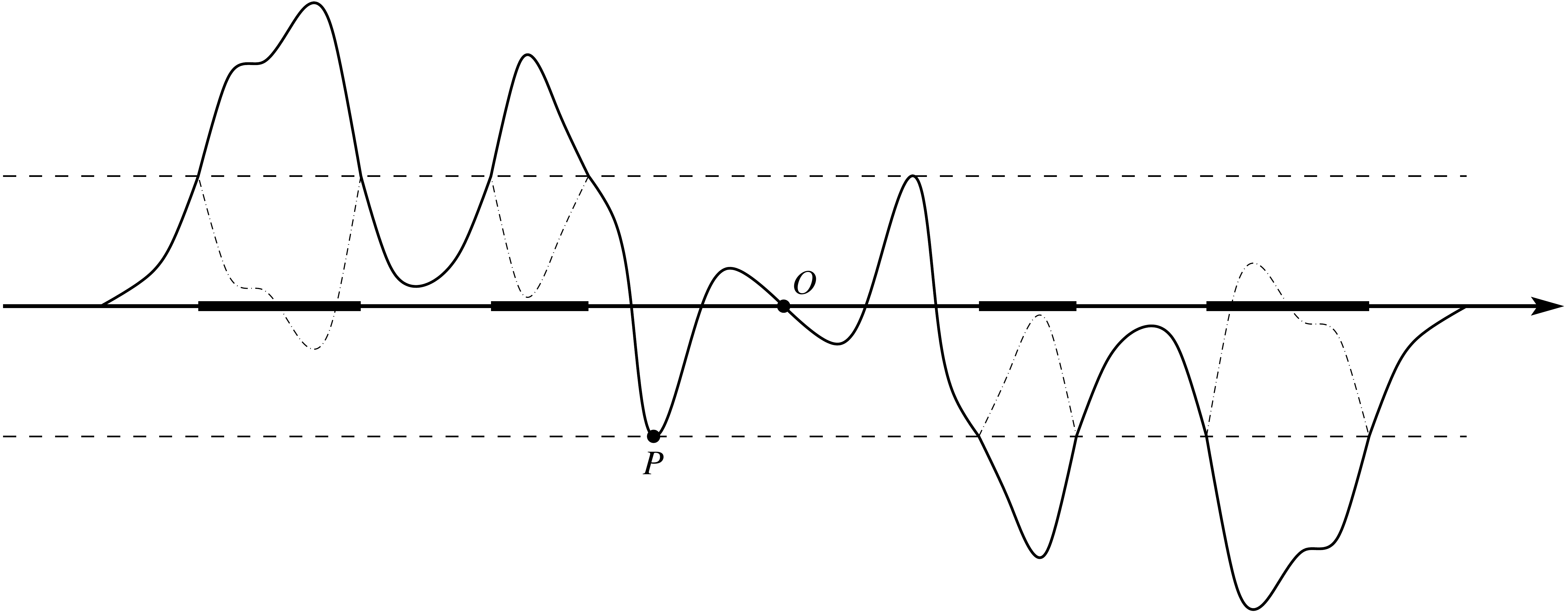}
	\caption{\sl The geometry involved in the proof of Lemma~\ref{SUSO}: detection of cancellations via isometric regions. In particular, the sets~$F_-$ and~$F_+$ in~\eqref{ieuryurfbsfbsjffjfghfweeu09998ye73}
	are given by the bold intervals along the
horizontal axis.}
	\label{DIFI2}
\end{figure}

We notice that, in light of~\eqref{SMG-2}, \eqref{SMG-8} and~\eqref{UYP},
\begin{equation*} F_-\subseteq(0,+\infty)\quad{\mbox{ and }}\quad
F_+\subseteq (-\infty,0).\end{equation*}
Furthermore, by the odd symmetry of~$u$ in~\eqref{SMG-1}, we have that
\begin{equation}\label{CHECK0}
{\mbox{$x\in F_-$ if and only if~$-x\in F_+$.}}\end{equation}

The strategy now is based on the detection of suitable cancellations via isometric regions
that correspond to the sets~$F_-$ and~$F_+$.
The gist is to get rid of the contributions arising from the complement of~$E$ below the line~$\{y=u(p)\}$. For this,
one has to use suitable transformations inherited from the geometry of the problem, such as
an odd reflection through the origin, a vertical translation of magnitude~$2u(p)$ and the reflection
along the line~$\{y=u(p)\}$. The use of these transformations will aim, on the one hand, at
detecting isometric regions in~$E$ which cancel the ones in the complement of~$E$ in the integral
computations, and, on the other hand, at maintaining a favorable control of the distance from the point~$P$, since
this quantity appears at the denominator of the integrands involved.

To employ this strategy, we use the notation~$X=(x,y)$
and the transformation
$$\widetilde{X}=(\widetilde x,\widetilde y)
=\widetilde{X}(X)
:=\big(-x,2u(p)-y\big)$$
and we claim that
\begin{equation}\label{CHECK1}
{\mbox{if~$x\ge0$, then~$|X-P|\ge|\widetilde X-P|$.}}
\end{equation}
To check this, we observe that when~$x\ge0$ it follows that~$xp\le0$ and therefore
\begin{eqnarray*}&&
|\widetilde X-P|^2-|X-P|^2\\&=&\Big| \big(-x,2u(p)-y\big)-\big(p,u(p)\big)\Big|^2
-\Big| (x,y)-\big(p,u(p)\big)\Big|^2\\&=&
\Big| \big(-x-p,u(p)-y\big)\Big|^2
-\Big| \big(x-p,y-u(p)\big)\Big|^2\\
&=&(-x-p)^2-(x-p)^2\\&=& 4xp\\&\le&0,
\end{eqnarray*}
which establishes~\eqref{CHECK1}.

Besides, we define
\begin{eqnarray*}
&&G:=\Big\{ X=(x,y)\in\R^2\setminus E {\mbox{ s.t. $x\in F_-$ and $y<2u(p)-u(x)$}}\Big\}
\end{eqnarray*}
and we claim that
\begin{equation}\label{CHECK2}
{\mbox{if~$X\in G$, then~$\widetilde X\in E\cap\{\widetilde x\in F_+\}\cap\{-u(\widetilde x)<
\widetilde y<2u(p)+u(\widetilde x)\}$.}}
\end{equation}
Indeed, if~$X\in G$, by the odd symmetry in~\eqref{SMG-1} we have that~$y>u(x)=-u(-x)$ and accordingly
\begin{equation}\label{CHECK05}
\widetilde y=2u(p)-y<2u(p)+u(-x)=2u(p)+u(\widetilde x).\end{equation}
In particuar, $\widetilde y<u(\widetilde x)$.
{F}rom this, \eqref{CHECK0} and~\eqref{CHECK05} we obtain that
\begin{equation}\label{CHECK7}
\widetilde X\in E\cap\{\widetilde x\in F_+\}\cap\{\widetilde y<2u(p)+u(\widetilde x)\}.\end{equation}
Additionally, 
$$ \widetilde y=2u(p)-y>2u(p)-\big(2u(p)-u(x) \big)=u(x)=u(-\widetilde x)=-u(\widetilde x).$$
Combining this with~\eqref{CHECK7} we obtain~\eqref{CHECK2}, as desired.

Now, in light of~\eqref{CHECK1} and~\eqref{CHECK2}, and observing that~$d\widetilde X=dX$,
\begin{equation}\label{CONT1}
\begin{split}
&\int_{\R^2\cap\{x\in F_-\}\cap\{y<2u(p)-u(x)\}}\frac{\chi_{\R^2\setminus E}(X)}{|X-P|^{2+s}}\,dX\\&\qquad=
\int_{G}\frac{dX}{|X-P|^{2+s}}\le
\int_{G}\frac{dX}{|\widetilde X(X)-P|^{2+s}}\\&\qquad\le
\int_{E\cap\{\widetilde x\in F_+\}\cap\{-u(\widetilde x)<
\widetilde y<2u(p)+u(\widetilde x)\}}\frac{d\widetilde X}{|\widetilde X-P|^{2+s}}\\&\qquad=\int_{
\R^2\cap\{\widetilde x\in F_+\}\cap\{-u(\widetilde x)<
\widetilde y<2u(p)+u(\widetilde x)\}}\frac{\chi_E(\widetilde X)}{|\widetilde X-P|^{2+s}}\,d\widetilde X.
\end{split}\end{equation}

Now we use the even reflection across the line~$\{y=u(p)\}$, that is we define
\begin{equation}\label{CONT10}X_*=(x_*,y_*)=X_*(X):=\big(x,2u(p)-y\big).\end{equation}
We notice that
\begin{equation}\label{CONT9} |X_*-P|=\Big| \big(x,2u(p)-y\big)-\big(p,u(p)\big)\Big|
=\Big| \big(x-p,u(p)-y\big)\Big|=|X-P|.\end{equation}
As a result,
\begin{eqnarray*}&&
\int_{\R^2\cap\{x\in F_-\}\cap\{y>2u(p)-u(x)\}}\frac{\chi_{\R^2\setminus E}(X)}{|X-P|^{2+s}}\,dX\le
\int_{\R^2\cap\{x\in F_-\}\cap\{y>2u(p)-u(x)\}}\frac{dX}{|X-P|^{2+s}}
\\&&\quad=
\int_{\R^2\cap\{x\in F_-\}\cap\{y>2u(p)-u(x)\}}\frac{dX}{|X_*(X)-P|^{2+s}}=
\int_{\R^2\cap\{x\in F_-\}\cap\{y_*<u(x_*)\}}\frac{dX_*}{|X_*-P|^{2+s}}\\&&\qquad=
\int_{\R^2\cap\{x\in F_-\}}\frac{\chi_E(X_*)}{|X_*-P|^{2+s}}\,dX_*.
\end{eqnarray*}
{F}rom this and~\eqref{CONT1}, changing the names of the integration variables, we arrive at
\begin{eqnarray*}&&
\int_{\R^2\cap\{x\in F_-\}}\frac{\chi_{\R^2\setminus E}(X)}{|X-P|^{2+s}}\,dX\\&&\qquad\le
\int_{\R^2\cap\{ x\in F_+\}\cap\{-u( x)< y<2u(p)+u( x)\}}\frac{\chi_E( X)}{| X-P|^{2+s}}\,dX
+\int_{\R^2\cap\{x\in F_-\}}\frac{\chi_E(X)}{|X-P|^{2+s}}\,dX,
\end{eqnarray*}
that is
\begin{equation*}
\int_{\R^2\cap\{x\in F_-\}}\frac{\chi_{\R^2\setminus E}(X)-\chi_E(X)}{|X-P|^{2+s}}\,dX\le
\int_{\R^2\cap\{ x\in F_+\}\cap\{-u( x)< y<2u(p)+u( x)\}}\frac{\chi_E( X)}{| X-P|^{2+s}}\,dX.
\end{equation*}

Therefore,
\begin{equation}\label{CONT22}\begin{split}&
\int_{\R^2\cap\{x\in F_-\cup F_+\}}\frac{\chi_{\R^2\setminus E}(X)-\chi_E(X)}{|X-P|^{2+s}}\,dX\\ \le&
\int_{\R^2\cap\{x\in F_+\}}\frac{\chi_{\R^2\setminus E}(X)-\chi_E(X)}{|X-P|^{2+s}}\,dX
+
\int_{\R^2\cap\{ x\in F_+\}\cap\{-u( x)< y<2u(p)+u( x)\}}\frac{\chi_E( X)}{| X-P|^{2+s}}\,dX
\\=&\int_{\R^2\cap\{x\in F_+\}}\frac{\chi_{\R^2\setminus E}(X)}{|X-P|^{2+s}}\,dX
-
\int_{\R^2\cap\{ x\in F_+\}\cap\{y\in(-\infty,-u( x))\cup(2u(p)+u( x),+\infty)\}}\frac{\chi_E( X)}{| X-P|^{2+s}}\,dX.
\end{split}\end{equation}

Now, recalling the notation in~\eqref{CONT10}, we note that
\begin{equation}\label{CONT7}
{\mbox{if~$X\in(\R^2\setminus E)\cap\{x\in F_+\}$, then~$X_*\in E\cap\{x_*\in F_+\}\cap\{y_*\in(-\infty,-u( x_*))\}$.}}\end{equation}
Indeed, if~$X$ is as above then, since~$u(p)<0$,
\begin{equation} \label{CONT8}y_*=2u(p)-y<2u(p)-u(x)<-u(x)=-u(x_*).\end{equation}
Hence, since~$u(x_*)=u(x)>-u(p)>0$,
$$ y_*<-u(x_*)<u(x_*).
$$
By gathering this and~\eqref{CONT8} we obtain the desired result in~\eqref{CONT7}.

Thus, by~\eqref{CONT9} and~\eqref{CONT7} we deduce that
\begin{eqnarray*}
&&\int_{\R^2\cap\{x\in F_+\}}\frac{\chi_{\R^2\setminus E}(X)}{|X-P|^{2+s}}\,dX
=\int_{(\R^2\setminus E)\cap\{x\in F_+\}}\frac{dX}{|X_*(X)-P|^{2+s}}\\&&\qquad\le
\int_{E\cap\{x_*\in F_+\}\cap\{y_*\in(-\infty,-u( x_*))\}}\frac{dX_*}{|X_*-P|^{2+s}}.
\end{eqnarray*}
Note that we have used here the notation~$X_*(X)$ to emphasize the dependence of~$X_*$ on the original coordinate~$X$ and use a change of variable in the integral calculation.
We write the above inequality in the form
\begin{eqnarray*}
\int_{\R^2\cap\{x\in F_+\}}\frac{\chi_{\R^2\setminus E}(X)}{|X-P|^{2+s}}\,dX
\le
\int_{\R^2\cap\{x\in F_+\}\cap\{y\in(-\infty,-u( x))\}}\frac{\chi_E(X)}{|X-P|^{2+s}}\,dX.
\end{eqnarray*}
Comparing with~\eqref{CONT22} we thereby deduce that
\begin{equation}\label{CONT23}
\int_{\R^2\cap\{x\in F_-\cup F_+\}}\frac{\chi_{\R^2\setminus E}(X)-\chi_E(X)}{|X-P|^{2+s}}\,dX\le0.\end{equation}

We also remark that
$$ \R\setminus(F_-\cup F_+)=\big\{x\in\R {\mbox{ s.t. }} |u(x)|\le-u(p)\big\}.$$
{F}rom this observation, \eqref{S:10} and~\eqref{CONT23} we infer that
\begin{equation}\label{S:10bi}
0\le\int_{\R^2\cap\{ |u(x)|\le-u(p)\}}\frac{\chi_{\R^2\setminus E}(X)-\chi_{E}(X)}{|X-P|^{2+s}}\,dX.\end{equation}

It is now useful to remark that
\begin{equation*}
{\mbox{if~$|u(x)|\le-u(p)$ and~$y<u(p)$, then~$X=(x,y)\in E$,}}
\end{equation*}
since in this setting~$u(x)\ge u(p)>y$.

Consequently, the inequality in~\eqref{S:10bi} gives that
\begin{eqnarray*}
0&\le&\int_{\R^2\cap\{ |u(x)|\le-u(p)\}\cap \{ y > u(p) \}}\frac{\chi_{\R^2\setminus E}(X)-\chi_{E}(X)}{|X-P|^{2+s}}\,dX\\&&
\qquad
+\int_{\R^2\cap\{ |u(x)|\le-u(p)\}\cap \{ y<u(p) \}}\frac{\chi_{\R^2\setminus E}(X)-\chi_{E}(X)}{|X-P|^{2+s}}\,dX\\
&=&
\int_{\R^2\cap\{ |u(x)|\le-u(p)\}\cap \{ y > u(p) \}}\frac{dX}{|X-P|^{2+s}}-2
\int_{E\cap\{ |u(x)|\le-u(p)\}\cap \{ y > u(p) \}}\frac{dX}{|X-P|^{2+s}}
\\&&\qquad-
\int_{\R^2\cap\{ |u(x)|\le-u(p)\}\cap \{ y<u(p) \}}\frac{dX}{|X-P|^{2+s}}.
\end{eqnarray*}
Hence, recalling~\eqref{CONT9} and noticing that, if~$y\ge u(p)$ then~$u(y_*)\le u(p)$, we find that
\begin{eqnarray*}
0&\le&
\int_{\R^2\cap\{ |u(x_*)|\le-u(p)\}\cap \{ y_* < u(p) \}}\frac{dX_*}{|X_*-P|^{2+s}}-2
\int_{E\cap\{ |u(x)|\le-u(p)\}\cap \{ y > u(p) \}}\frac{dX}{|X-P|^{2+s}}
\\&&\qquad-
\int_{\R^2\cap\{ |u(x)|\le-u(p)\}\cap \{ y<u(p) \}}\frac{dX}{|X-P|^{2+s}}\\&=&-2
\int_{E\cap\{ |u(x)|\le-u(p)\}\cap \{ y > u(p) \}}\frac{dX}{|X-P|^{2+s}}
.
\end{eqnarray*}
This entails that the set
\begin{equation}\label{PoryrgfgfvhgreiYYYY}
{\mbox{$E\cap\{ |u(x)|\le-u(p)\}\cap \{ y > u(p) \}$ is of null measure.}}\end{equation}

On the other hand, we set
$$ \eta:=\min\left\{\frac{|u(p)|}{1+8\|u\|_{C^1(-d/2,d/2)}}, \frac{d}2\right\}$$
and
we claim that
\begin{equation}\label{PoryrgfgfvhgreiYYYY2}
\left[ -\eta, \eta\right]\times
\left[\frac{u(p)}2, \frac{u(p)}4\right] \subseteq E\cap\{ |u(x)|\le-u(p)\}\cap \{ y > u(p) \}.
\end{equation}
Indeed, if~$(x,y)\in\left[ -\eta, \eta\right]\times
\left[\frac{u(p)}2, \frac{u(p)}4\right]$, then
\begin{eqnarray*}&&
y-u(x)<\frac{u(p)}4-u(x) +u(0)\le \frac{u(p)}4+\|u\|_{C^1(-d/2,d/2)} |x|\\&&\qquad
\le \frac{u(p)}4+ \frac{|u(p)|\,\|u\|_{C^1(-d/2,d/2)}}{1+8\|u\|_{C^1(-d/2,d/2)}}
\le \frac{u(p)}4+ \frac{|u(p)| }{8}\le  \frac{u(p)}8<0,
\end{eqnarray*}
which gives that~$(x,y)\in E$.

Moreover,
\begin{eqnarray*}&&
|u(x)|=|u(x)-u(0)|\le \|u\|_{C^1(-d/2,d/2)} |x|\le \frac{|u(p)|\,\|u\|_{C^1(-d/2,d/2)}}{1+8\|u\|_{C^1(-d/2,d/2)}}\le |u(p)|
\end{eqnarray*}
and~$y\ge \frac{u(p)}2 >u(p)$. These considerations prove~\eqref{PoryrgfgfvhgreiYYYY2}.

{F}rom~\eqref{PoryrgfgfvhgreiYYYY} and~\eqref{PoryrgfgfvhgreiYYYY2} we obtain the desired contradiction.
The proof of~\eqref{SMG-4} is thereby complete.\end{proof}

We are now in position of completing the proof of the antisymmetric maximum principle in
Theorem~\ref{SMG-T}. \begin{proof}[Proof of Theorem~\ref{SMG-T}]
We observe that the claim in~\eqref{SMG-3} about the odd
symmetry of~$u$ is a consequence of the antisymmetric property
of~$u$ outside~$(-d,d)$, as given by~\eqref{SMG-1TH},
and~\cite[Lemma~A.1]{MR3596708}.

Thus, to complete the proof of Theorem~\ref{SMG-T} it remains to establish~\eqref{920irpkfeitk-320jtjoyj-01}
(with this, the claim in~\eqref{920irpkfeitk-320jtjoyj-02} would then follow from the odd symmetry of~$u$
demonstrated in~\eqref{SMG-3}).

Hence, we can assume that there exists~$p_\star\in[-d,0]$ such that
\begin{equation}\label{ECX58i2p-3}
u(p_\star)<0,
\end{equation}
otherwise~\eqref{920irpkfeitk-320jtjoyj-01} would be automatically satisfied.

We now aim at applying Lemma~\ref{SUSO}, which would entail~\eqref{920irpkfeitk-320jtjoyj-01}
and thus end the proof of Theorem~\ref{SMG-T}. For this, we need to check that all the hypotheses of Lemma~\ref{SUSO} are fulfilled.

{F}rom~\eqref{SMG-3}, we have that condition~\eqref{SMG-1} is satisfied.
Moreover, condition~\eqref{SMG-2} holds true, due to~\eqref{SMG-2TH}.

We recall also that~$u$ is uniformly continuous in~$(-d,d)$,
owing to~\cite[Theorem~1.1]{MR3516886}. 
As a consequence, we can redefine~$u$ at the extrema of~$(-d,d)$ by setting
\begin{equation*} u(-d):=\lim_{x\searrow-d}u(x)\quad{\mbox{and}}\quad
u(d):=\lim_{x\nearrow d}u(x)\end{equation*}
and we have that
\begin{equation}\label{SMG-7}
u\in C([-d,d]).\end{equation}
Furthermore, we recall that
\begin{equation}\label{0ujo9u0t4g0DZpfejhf76nh}
u\in C^\infty(-d,d),\end{equation} thanks to~\cite{MR3090533, MR3331523}.
This and~\eqref{SMG-7} give that the regularity assumptions on~$u$ taken in Lemma~\ref{SUSO}
are satisfied.

Also, in light of~\eqref{0ujo9u0t4g0DZpfejhf76nh},
we can write the Euler-Lagrange equation for nonlocal minimal surfaces (see~\cite[Theorem~5.1]{MR2675483})
in a pointwise sense in~$(-d,d)$. In particular, condition~\eqref{S:10} is fulfilled as well.

\begin{figure}
		\centering
		\includegraphics[width=.4\linewidth]{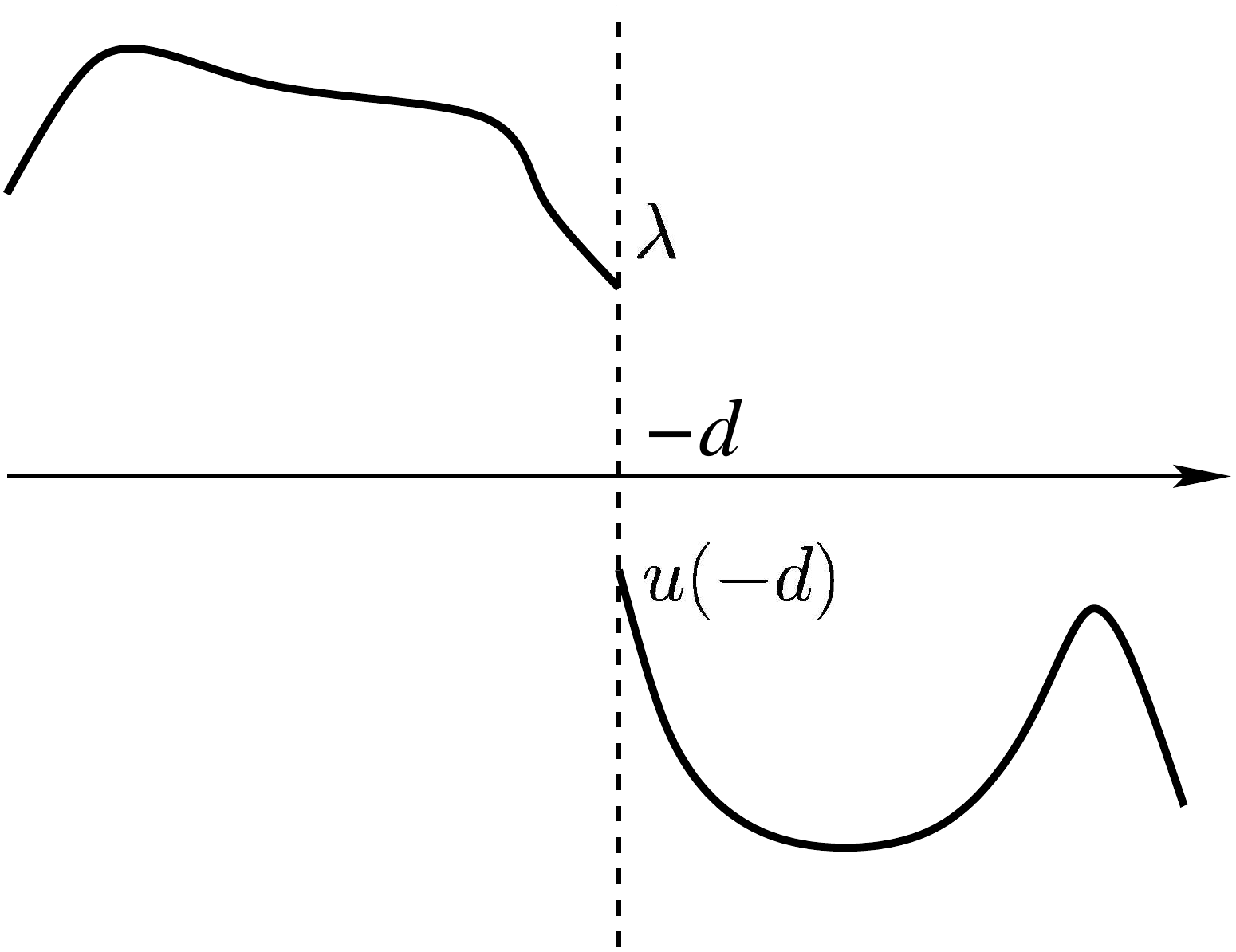}
	\caption{\sl The geometry involved in the proof of Theorem~\ref{SMG-T}: detachment from the boundary.}
			\label{DIFI}
\end{figure}

Consequently, to apply Lemma~\ref{SUSO}, it remains to check condition~\eqref{SECCA}.
For this, suppose by contradiction that~$u(-d)=\min_{[-d,0]}u$.
This, together with the assumptions on~$u$ in~\eqref{SMG-1TH} and~\eqref{SMG-2TH},
and recalling also~\eqref{ECX58i2p-3},
gives that
\begin{equation}\label{SMG-87} \lambda:=\lim_{x\nearrow-d}u(x)\ge 0>u(p_\star)\ge\min_{[-d,0]}u=u(-d),\end{equation}
see Figure~\ref{DIFI}.

Hence, by~\cite[Corollary~1.3(ii)]{MR4104542} (see also~\cite{MR3532394} for related results),
we infer that there exist~$\mu>0$ 
and a continuous function~$v$ in~$\big(u(-d)-\mu,u(-d)+\mu\big)$
such that~$v$ is the inverse function of~$u$ (in the above domain of definition, with the understanding that~$v=-d$
along the jump discontinuity of~$u$ detected in~\eqref{SMG-87}).

Now we pick a sequence of points~$q_k\in(-d,0)$ with~$q_k\searrow-d$ as~$k\to+\infty$;
in particular, we can assume that~$u(q_k)\in\big(u(-d)-\mu,u(-d)+\mu\big)$
and we thus obtain that~$v(u(q_k))=q_k$. 
We have that
\begin{equation}\label{SMG-86}
u(q_k)< u(-d),
\end{equation}
otherwise, since~$v(y)=-d$ for all~$y\in\big[u(-d),\min\{\lambda,u(-d)+\mu\}\big]$, we would have that~$-d=v(u(q_k))=q_k>-d$,
which is a contradiction.

Having proved~\eqref{SMG-86}, we find that
$$u(q_k)< u(-d)=u(p)=\min_{[-d,0]}u\le u(q_k)$$
and this is a contradiction that shows the validity of~\eqref{SECCA}.

We can thereby apply Lemma~\ref{SUSO} to complete the proof of Theorem~\ref{SMG-T}.\end{proof}

\section{Some comments about antisymmetric functions}\label{SNEW}

Maximum principles in the presence of an odd symmetry have been already treated in the field
of integro-differential equations (see for instance~\cite{MR3453602, MR4030266, MR4108219, MR4308250, JACK}).
The main idea in this setting is that the operator can be rewritten as a different
integro-differential operator acting only on functions defined in the halfspace,
but still with a positive kernel (under certain assumptions on the original operator).
For instance, the fractional Laplacian in~$\R$ acting on odd functions~$u$
can be rewritten as an operator~$L_su+cu$, where~$c$ is a nonnegative function
and~$L_s$ is an integro-differential operator of the form
\begin{equation}\label{ESNEW}  L_s u(x) =\int_0^{+\infty}\big(u(x)-u(y)\big) K_s(x,y)\,dy \end{equation}
defined for all~$x >0$,  being~$K_s(x,y)$ a nonnegative kernel.
Thus, exploiting this observation, the usual proof for the maximum principle in~$\R$ carries through also for this operator.

A natural question is whether or not a trick of this type, exploiting the antisymmetry of the functions,
can work in this setting for antisymmetric nonlocal minimal graphs. 

In a sense, the answer to this question is not completely clear:
on the one hand, there are ways to rewrite the nonlocal mean curvature operator
in the odd symmetric setting,
but either the kernel obtained has not necessarily the right sign, or the structure of the
operator is not immediately
apt for a straightforward proof of the maximum principle.

In our opinion, the difference between the fractional Laplace framework and that of
nonlocal minimal surfaces
with respect to this point may be not merely algebraic in nature and instead reveals an
interesting phenomenon
due to the nonlinear and geometric structures of the problem (thus requiring a somewhat 
careful detection
of the appropriate isometric regions in Section~\ref{SEC-PF2}).

\begin{figure}
		\centering
		\includegraphics[width=.9\linewidth]{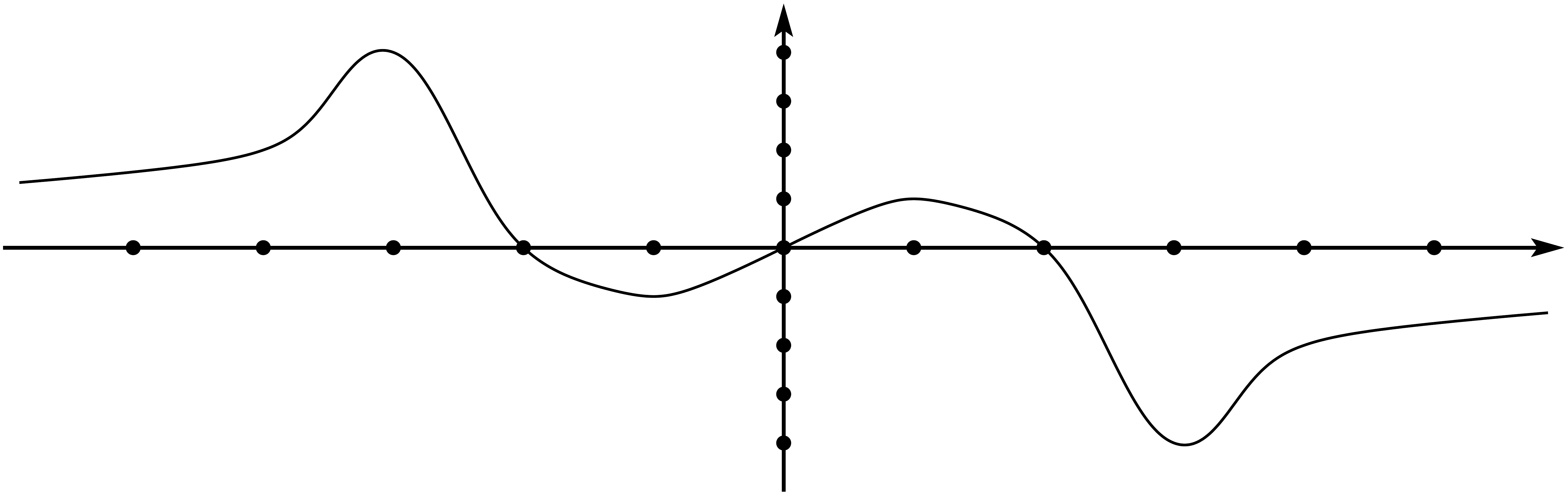}
	\caption{\sl Difficulties arising in a straightforward proof of the antisymmetric maximum principle.}
	\label{3DIFI2000BIS}
\end{figure}

Let us explicitly point out some of the natural computations that one can perform to encode 
odd symmetry properties into the nonlocal mean curvature operator. First of all, at the level 
of a set~$E$, the antisymmetric property considered in this paper
reads that~$X\in E$ if and only if~$-X\in\R^2\setminus E$. Accordingly,
$$ \chi_{\R^2\setminus E}(-X)-\chi_{E}(-X)=
-\big( \chi_{\R^2\setminus E}(X)-\chi_{E}(X)\big),$$
leading to
\begin{eqnarray*}
&&\int_{\R^2}\frac{\chi_{\R^2\setminus E}(X)-\chi_{E}(X)}{|X-P|^{2+s}}\,dX\\&=&
\int_{(0,+\infty)\times\R}\frac{\chi_{\R^2\setminus E}(X)-\chi_{E}(X)}{|X-P|^{2+s}}\,dX+
\int_{(-\infty,0)\times\R}\frac{\chi_{\R^2\setminus E}(X)-\chi_{E}(X)}{|X-P|^{2+s}}\,dX\\&=&
\int_{(0,+\infty)\times\R}\frac{\chi_{\R^2\setminus E}(X)-\chi_{E}(X)}{|X-P|^{2+s}}\,dX+
\int_{(0,+\infty)\times\R}\frac{\chi_{\R^2\setminus E}(-X)-\chi_{E}(-X)}{|X+P|^{2+s}}\,dX\\&=&
\int_{(0,+\infty)\times\R}\frac{\chi_{\R^2\setminus E}(X)-\chi_{E}(X)}{|X-P|^{2+s}}\,dX-
\int_{(0,+\infty)\times\R}\frac{\chi_{\R^2\setminus E}(X)-\chi_{E}(X)}{|X+P|^{2+s}}\,dX\\&=&
\int_{(0,+\infty)\times\R}\big( \chi_{\R^2\setminus E}(X)-\chi_{E}(X)\big)\,K(X,P)\,dX,
\end{eqnarray*}
where
$$ K(X,P):=\frac1{|X-P|^{2+s}}-\frac1{|X+P|^{2+s}}.$$
While this calculation reduces the nonlocal mean curvature to an integral in the halfplane, the kernel~$K$ is not necessarily positive,
for instance, if~$P=(1,1)$ and~$X=(1,-2)$, we have that
$$K(X,P)=\frac{1}{3^{2+s}}-\frac{1}{5^{\frac{2+s}2}}<0 ,$$ see Figure~\ref{3DIFI2000BIS}.

Another possible approach towards the rewriting of the nonlocal mean curvature operator
consists in focusing on sets possessing a graph structure (say, in the plane, for simplicity).
In this case, by equation~(49) in~\cite{MR3331523},
there exists a bounded, smooth, globally Lipschitz, strictly increasing and odd function~$F$ such that,
up to a real multiplicative constant, the nonlocal mean curvature at~$(x,u(x))$ becomes
\begin{eqnarray*}
\int_\R F\left(\frac{u(x)-u(x+w)}{|w|}\right)\,\frac{dw}{|w|^{1+s}}=
\int_\R F\left(\frac{u(x)-u(y)}{|x-y|}\right)\,\frac{dy}{|x-y|^{1+s}}.
\end{eqnarray*}
Thus, if~$u$ is odd symmetric, we obtain
\begin{eqnarray*}
&&\int_0^{+\infty} F\left(\frac{u(x)-u(y)}{|x-y|}\right)\,\frac{dy}{|x-y|^{1+s}}
+\int_0^{+\infty} F\left(\frac{u(x)-u(-y)}{|x+y|}\right)\,\frac{dy}{|x+y|^{1+s}}\\&=&
\int_0^{+\infty} F\left(\frac{u(x)-u(y)}{|x-y|}\right)\,\frac{dy}{|x-y|^{1+s}}
-\int_0^{+\infty} F\left(\frac{u(x)-u(y) -2u(x)}{|x+y|}\right)\,\frac{dy}{|x+y|^{1+s}},
\end{eqnarray*}
which, due to the nonlinearity~$F$, does not immediately lead to a nice equation as in~\eqref{ESNEW}.

Yet, one could try to modify this argument by using Lagrange's Mean Value Theorem: in this way,
we can write the nonlocal mean curvature operator as
\begin{eqnarray*}&&
\int_0^{+\infty} F\left(\frac{u(x)-u(x+w)}{|w|}\right)\,\frac{dw}{|w|^{1+s}}
+\int_0^{+\infty} F\left(\frac{u(x)-u(x-w)}{|w|}\right)\,\frac{dw}{|w|^{1+s}}\\&=&
\int_0^{+\infty} F\left(\frac{u(x)-u(x+w)}{|w|}\right)\,\frac{dw}{|w|^{1+s}}-\int_0^{+\infty} F\left(\frac{u(x-w)-u(x)}{|w|}\right)\,\frac{dw}{|w|^{1+s}}\\&=&
\int_0^{+\infty}\frac{F'(\xi_u(x,w))}{|w|^{2+s}}\,
\Big( 2u(x)-u(x+w)-u(x-w)\Big)\,dw,
\end{eqnarray*}
even without the odd symmetry of~$u$. In this case, the kernel~$\frac{F'(\xi_u(x,w))}{|w|^{2+s}}$ is positive,
since~$F$ is strictly increasing, and bounded from above by the kernel of the fractional Laplacian, since $F$ is globally Lipschitz.

Nonetheless, if~$u$ is odd symmetric and, say, negative in~$(2,+\infty)$ and presents a positive maximum in~$[0,2]$,
it is not necessarily true that~$2u(x)-u(x+w)-u(x-w)\ge0$. As an example, one can consider an odd function as above such that~$u(1)=1$, $u(3)=-4$
and~$u(4)=-2$: in this way, if~$x:=1$ and~$w:=4$, it follows that~$2u(x)-u(x+w)-u(x-w)=
2u(1)-u(5)-u(-3)=2+1-4<0$, see Figure~\ref{3DIFI2000BIS}.

For these reasons, it is not clear to us that a straightforward symmetric trick may lead to 
an antisymmetric maximum principle as in Theorem~\ref{SMG-T} without taking advantage of the cancellations
provided by the isometric regions detected in Section~\ref{SEC-PF2}.

In any case, we think that a proof of the antisymmetric maximum principle based
solely on the detection of isometric regions, as presented in Section~\ref{SEC-PF2},
is interesting in itself, since it explicitly identifies the convenient cancellations
for such a result to hold.

\vfill


\begin{bibdiv}
\begin{biblist}

\bib{MR3331523}{article}{
   author={Barrios, Bego\~{n}a},
   author={Figalli, Alessio},
   author={Valdinoci, Enrico},
   title={Bootstrap regularity for integro-differential operators and its
   application to nonlocal minimal surfaces},
   journal={Ann. Sc. Norm. Super. Pisa Cl. Sci. (5)},
   volume={13},
   date={2014},
   number={3},
   pages={609--639},
   issn={0391-173X},
   review={\MR{3331523}},
}

\bib{MR3982031}{article}{
   author={Borthagaray, Juan Pablo},
   author={Li, Wenbo},
   author={Nochetto, Ricardo H.},
   title={Finite element discretizations of nonlocal minimal graphs:
   convergence},
   journal={Nonlinear Anal.},
   volume={189},
   date={2019},
   pages={111566, 31},
   issn={0362-546X},
   review={\MR{3982031}},
   doi={10.1016/j.na.2019.06.025},
}

\bib{MR4294645}{article}{
   author={Borthagaray, Juan Pablo},
   author={Li, Wenbo},
   author={Nochetto, Ricardo H.},
   title={Finite element algorithms for nonlocal minimal graphs},
   journal={Math. Eng.},
   volume={4},
   date={2022},
   number={2},
   pages={Paper No. 016, 29},
   review={\MR{4294645}},
   doi={10.3934/mine.2022016},
}

\bib{MR4184583}{article}{
   author={Bucur, Claudia},
   author={Dipierro, Serena},
   author={Lombardini, Luca},
   author={Valdinoci, Enrico},
   title={Minimisers of a fractional seminorm and nonlocal minimal surfaces},
   journal={Interfaces Free Bound.},
   volume={22},
   date={2020},
   number={4},
   pages={465--504},
   issn={1463-9963},
   review={\MR{4184583}},
   doi={10.4171/ifb/447},
}

\bib{MR3926519}{article}{
   author={Bucur, Claudia},
   author={Lombardini, Luca},
   author={Valdinoci, Enrico},
   title={Complete stickiness of nonlocal minimal surfaces for small values
   of the fractional parameter},
   journal={Ann. Inst. H. Poincar\'{e} Anal. Non Lin\'{e}aire},
   volume={36},
   date={2019},
   number={3},
   pages={655--703},
   issn={0294-1449},
   review={\MR{3926519}},
   doi={10.1016/j.anihpc.2018.08.003},
}

\bib{MR4116635}{article}{
   author={Cabr\'{e}, Xavier},
   author={Cinti, Eleonora},
   author={Serra, Joaquim},
   title={Stable $s$-minimal cones in $\Bbb{R}^3$ are flat for $s\sim 1$},
   journal={J. Reine Angew. Math.},
   volume={764},
   date={2020},
   pages={157--180},
   issn={0075-4102},
   review={\MR{4116635}},
   doi={10.1515/crelle-2019-0005},
}

\bib{MR3934589}{article}{
   author={Cabr\'{e}, Xavier},
   author={Cozzi, Matteo},
   title={A gradient estimate for nonlocal minimal graphs},
   journal={Duke Math. J.},
   volume={168},
   date={2019},
   number={5},
   pages={775--848},
   issn={0012-7094},
   review={\MR{3934589}},
   doi={10.1215/00127094-2018-0052},
}

\bib{MR3744919}{article}{
   author={Cabr\'{e}, Xavier},
   author={Fall, Mouhamed Moustapha},
   author={Weth, Tobias},
   title={Delaunay hypersurfaces with constant nonlocal mean curvature},
   language={English, with English and French summaries},
   journal={J. Math. Pures Appl. (9)},
   volume={110},
   date={2018},
   pages={32--70},
   issn={0021-7824},
   review={\MR{3744919}},
   doi={10.1016/j.matpur.2017.07.005},
}

\bib{MR3881478}{article}{
   author={Cabr\'{e}, Xavier},
   author={Fall, Mouhamed Moustapha},
   author={Sol\`a-Morales, Joan},
   author={Weth, Tobias},
   title={Curves and surfaces with constant nonlocal mean curvature: meeting
   Alexandrov and Delaunay},
   journal={J. Reine Angew. Math.},
   volume={745},
   date={2018},
   pages={253--280},
   issn={0075-4102},
   review={\MR{3881478}},
   doi={10.1515/crelle-2015-0117},
}

\bib{MR3532394}{article}{
   author={Caffarelli, Luis},
   author={De Silva, Daniela},
   author={Savin, Ovidiu},
   title={Obstacle-type problems for minimal surfaces},
   journal={Comm. Partial Differential Equations},
   volume={41},
   date={2016},
   number={8},
   pages={1303--1323},
   issn={0360-5302},
   review={\MR{3532394}},
   doi={10.1080/03605302.2016.1192646},
}

\bib{MR2675483}{article}{
   author={Caffarelli, Luis},
   author={Roquejoffre, Jean-Michel},
   author={Savin, Ovidiu},
   title={Nonlocal minimal surfaces},
   journal={Comm. Pure Appl. Math.},
   volume={63},
   date={2010},
   number={9},
   pages={1111--1144},
   issn={0010-3640},
   review={\MR{2675483}},
   doi={10.1002/cpa.20331},
}

\bib{MR2564467}{article}{
   author={Caffarelli, Luis A.},
   author={Souganidis, Panagiotis E.},
   title={Convergence of nonlocal threshold dynamics approximations to front
   propagation},
   journal={Arch. Ration. Mech. Anal.},
   volume={195},
   date={2010},
   number={1},
   pages={1--23},
   issn={0003-9527},
   review={\MR{2564467}},
   doi={10.1007/s00205-008-0181-x},
}

\bib{MR3107529}{article}{
   author={Caffarelli, Luis},
   author={Valdinoci, Enrico},
   title={Regularity properties of nonlocal minimal surfaces via limiting
   arguments},
   journal={Adv. Math.},
   volume={248},
   date={2013},
   pages={843--871},
   issn={0001-8708},
   review={\MR{3107529}},
   doi={10.1016/j.aim.2013.08.007},
}

\bib{MR3156889}{article}{
   author={Chambolle, Antonin},
   author={Morini, Massimiliano},
   author={Ponsiglione, Marcello},
   title={Minimizing movements and level set approaches to nonlocal
   variational geometric flows},
   conference={
      title={Geometric partial differential equations},
   },
   book={
      series={CRM Series},
      volume={15},
      publisher={Ed. Norm., Pisa},
   },
   date={2013},
   pages={93--104},
   review={\MR{3156889}},
   doi={10.1007/978-88-7642-473-1\_4},
}

\bib{MR3713894}{article}{
   author={Chambolle, Antonin},
   author={Novaga, Matteo},
   author={Ruffini, Berardo},
   title={Some results on anisotropic fractional mean curvature flows},
   journal={Interfaces Free Bound.},
   volume={19},
   date={2017},
   number={3},
   pages={393--415},
   issn={1463-9963},
   review={\MR{3713894}},
   doi={10.4171/IFB/387},
}

\bib{MR3981295}{article}{
   author={Cinti, Eleonora},
   author={Serra, Joaquim},
   author={Valdinoci, Enrico},
   title={Quantitative flatness results and $BV$-estimates for stable
   nonlocal minimal surfaces},
   journal={J. Differential Geom.},
   volume={112},
   date={2019},
   number={3},
   pages={447--504},
   issn={0022-040X},
   review={\MR{3981295}},
   doi={10.4310/jdg/1563242471},
}

\bib{MR3778164}{article}{
   author={Cinti, Eleonora},
   author={Sinestrari, Carlo},
   author={Valdinoci, Enrico},
   title={Neckpinch singularities in fractional mean curvature flows},
   journal={Proc. Amer. Math. Soc.},
   volume={146},
   date={2018},
   number={6},
   pages={2637--2646},
   issn={0002-9939},
   review={\MR{3778164}},
   doi={10.1090/proc/14002},
}

\bib{MR3836150}{article}{
   author={Ciraolo, Giulio},
   author={Figalli, Alessio},
   author={Maggi, Francesco},
   author={Novaga, Matteo},
   title={Rigidity and sharp stability estimates for hypersurfaces with
   constant and almost-constant nonlocal mean curvature},
   journal={J. Reine Angew. Math.},
   volume={741},
   date={2018},
   pages={275--294},
   issn={0075-4102},
   review={\MR{3836150}},
   doi={10.1515/crelle-2015-0088},
}

\bib{MR4108219}{article}{
   author={Cozzi, Matteo},
   author={D\'{a}vila, Juan},
   author={del Pino, Manuel},
   title={Long-time asymptotics for evolutionary crystal dislocation models},
   journal={Adv. Math.},
   volume={371},
   date={2020},
   pages={107242, 109},
   issn={0001-8708},
   review={\MR{4108219}},
   doi={10.1016/j.aim.2020.107242},
}

\bib{MR4291290}{article}{
   author={Cozzi, Matteo},
   author={Farina, Alberto},
   author={Lombardini, Luca},
   title={Bernstein-Moser-type results for nonlocal minimal graphs},
   journal={Comm. Anal. Geom.},
   volume={29},
   date={2021},
   number={4},
   pages={761--777},
   issn={1019-8385},
   review={\MR{4291290}},
   doi={10.4310/CAG.2021.v29.n4.a1},
}

\bib{MR4279395}{article}{
   author={Cozzi, Matteo},
   author={Lombardini, Luca},
   title={On nonlocal minimal graphs},
   journal={Calc. Var. Partial Differential Equations},
   volume={60},
   date={2021},
   number={4},
   pages={Paper No. 136, 72},
   issn={0944-2669},
   review={\MR{4279395}},
   doi={10.1007/s00526-021-02002-9},
}

\bib{MR3485130}{article}{
   author={D\'{a}vila, Juan},
   author={del Pino, Manuel},
   author={Dipierro, Serena},
   author={Valdinoci, Enrico},
   title={Nonlocal Delaunay surfaces},
   journal={Nonlinear Anal.},
   volume={137},
   date={2016},
   pages={357--380},
   issn={0362-546X},
   review={\MR{3485130}},
   doi={10.1016/j.na.2015.10.009},
}

\bib{MR3798717}{article}{
   author={D\'{a}vila, Juan},
   author={del Pino, Manuel},
   author={Wei, Juncheng},
   title={Nonlocal $s$-minimal surfaces and Lawson cones},
   journal={J. Differential Geom.},
   volume={109},
   date={2018},
   number={1},
   pages={111--175},
   issn={0022-040X},
   review={\MR{3798717}},
   doi={10.4310/jdg/1525399218},
}

\bib{MR3412379}{article}{
   author={Di Castro, Agnese},
   author={Novaga, Matteo},
   author={Ruffini, Berardo},
   author={Valdinoci, Enrico},
   title={Nonlocal quantitative isoperimetric inequalities},
   journal={Calc. Var. Partial Differential Equations},
   volume={54},
   date={2015},
   number={3},
   pages={2421--2464},
   issn={0944-2669},
   review={\MR{3412379}},
   doi={10.1007/s00526-015-0870-x},
}

\bib{MR4096831}{article}{
   author={Dipierro, Serena},
   author={Dzhugan, Aleksandr},
   author={Forcillo, Nicol\`o},
   author={Valdinoci, Enrico},
   title={Enhanced boundary regularity of planar nonlocal minimal graphs and
   a butterfly effect},
   language={English, with English and Italian summaries},
   conference={
      title={Something about nonlinear problems},
   },
   book={
      series={Bruno Pini Math. Anal. Semin.},
      volume={11},
      publisher={Univ. Bologna, Alma Mater Stud., Bologna},
   },
   date={2020},
   pages={44--67},
   review={\MR{4096831}},
}

\bib{2020arXiv201000798D}{article}{
 author={Dipierro, Serena},
   author={Onoue, Fumihiko},
   author={Valdinoci, Enrico},
   title={(Dis)connectedness of nonlocal minimal surfaces in a cylinder and
   a stickiness property},
   journal={Proc. Amer. Math. Soc.},
   volume={150},
   date={2022},
   number={5},
   pages={2223--2237},
   issn={0002-9939},
   review={\MR{4392355}},
   doi={10.1090/proc/15796},
}

\bib{MR3516886}{article}{
   author={Dipierro, Serena},
   author={Savin, Ovidiu},
   author={Valdinoci, Enrico},
   title={Graph properties for nonlocal minimal surfaces},
   journal={Calc. Var. Partial Differential Equations},
   volume={55},
   date={2016},
   number={4},
   pages={Art. 86, 25},
   issn={0944-2669},
   review={\MR{3516886}},
   doi={10.1007/s00526-016-1020-9},
}

\bib{MR3596708}{article}{
   author={Dipierro, Serena},
   author={Savin, Ovidiu},
   author={Valdinoci, Enrico},
   title={Boundary behavior of nonlocal minimal surfaces},
   journal={J. Funct. Anal.},
   volume={272},
   date={2017},
   number={5},
   pages={1791--1851},
   issn={0022-1236},
   review={\MR{3596708}},
   doi={10.1016/j.jfa.2016.11.016},
}

\bib{MR4104542}{article}{
   author={Dipierro, Serena},
   author={Savin, Ovidiu},
   author={Valdinoci, Enrico},
   title={Nonlocal minimal graphs in the plane are generically sticky},
   journal={Comm. Math. Phys.},
   volume={376},
   date={2020},
   number={3},
   pages={2005--2063},
   issn={0010-3616},
   review={\MR{4104542}},
   doi={10.1007/s00220-020-03771-8},
}

\bib{MR4178752}{article}{
   author={Dipierro, Serena},
   author={Savin, Ovidiu},
   author={Valdinoci, Enrico},
   title={Boundary properties of fractional objects: flexibility of linear
   equations and rigidity of minimal graphs},
   journal={J. Reine Angew. Math.},
   volume={769},
   date={2020},
   pages={121--164},
   issn={0075-4102},
   review={\MR{4178752}},
   doi={10.1515/crelle-2019-0045},
}

\bib{JACK}{article}{
   author={Dipierro, Serena},
   author={Thompson, Jack},
   author={Valdinoci, Enrico},
title={On the Harnack inequality for antisymmetric $s$-harmonic functions},
journal={Preprint},
date={2022},
}

\bib{MR4050198}{article}{
   author={Farina, Alberto},
   author={Valdinoci, Enrico},
   title={Flatness results for nonlocal minimal cones and subgraphs},
   journal={Ann. Sc. Norm. Super. Pisa Cl. Sci. (5)},
   volume={19},
   date={2019},
   number={4},
   pages={1281--1301},
   issn={0391-173X},
   review={\MR{4050198}},
}

\bib{MR4308250}{article}{
   author={Felipe-Navarro, Juan-Carlos},
   title={Uniqueness for linear integro-differential equations in the real
   line and applications},
   journal={Calc. Var. Partial Differential Equations},
   volume={60},
   date={2021},
   number={6},
   pages={Paper No. 220, 25},
   issn={0944-2669},
   review={\MR{4308250}},
   doi={10.1007/s00526-021-02084-5},
}

\bib{MR4030266}{article}{
   author={Felipe-Navarro, Juan-Carlos},
   author={Sanz-Perela, Tom\'{a}s},
   title={Semilinear integro-differential equations, I: Odd solutions with
   respect to the Simons cone},
   journal={J. Funct. Anal.},
   volume={278},
   date={2020},
   number={2},
   pages={108309, 48},
   issn={0022-1236},
   review={\MR{4030266}},
   doi={10.1016/j.jfa.2019.108309},
}

\bib{MR3322379}{article}{
   author={Figalli, Alessio},
   author={Fusco, Nicola},
   author={Maggi, Francesco},
   author={Millot, Vincent},
   author={Morini, Massimiliano},
   title={Isoperimetry and stability properties of balls with respect to
   nonlocal energies},
   journal={Comm. Math. Phys.},
   volume={336},
   date={2015},
   number={1},
   pages={441--507},
   issn={0010-3616},
   review={\MR{3322379}},
   doi={10.1007/s00220-014-2244-1},
}

\bib{MR3680376}{article}{
   author={Figalli, Alessio},
   author={Valdinoci, Enrico},
   title={Regularity and Bernstein-type results for nonlocal minimal
   surfaces},
   journal={J. Reine Angew. Math.},
   volume={729},
   date={2017},
   pages={263--273},
   issn={0075-4102},
   review={\MR{3680376}},
   doi={10.1515/crelle-2015-0006},
}

\bib{MR2799577}{article}{
   author={Fusco, Nicola},
   author={Millot, Vincent},
   author={Morini, Massimiliano},
   title={A quantitative isoperimetric inequality for fractional perimeters},
   journal={J. Funct. Anal.},
   volume={261},
   date={2011},
   number={3},
   pages={697--715},
   issn={0022-1236},
   review={\MR{2799577}},
   doi={10.1016/j.jfa.2011.02.012},
}

\bib{MR2487027}{article}{
   author={Imbert, Cyril},
   title={Level set approach for fractional mean curvature flows},
   journal={Interfaces Free Bound.},
   volume={11},
   date={2009},
   number={1},
   pages={153--176},
   issn={1463-9963},
   review={\MR{2487027}},
   doi={10.4171/IFB/207},
}

\bib{MR3453602}{article}{
   author={Jarohs, Sven},
   author={Weth, Tobias},
   title={Symmetry via antisymmetric maximum principles in nonlocal problems
   of variable order},
   journal={Ann. Mat. Pura Appl. (4)},
   volume={195},
   date={2016},
   number={1},
   pages={273--291},
   issn={0373-3114},
   review={\MR{3453602}},
   doi={10.1007/s10231-014-0462-y},
}

\bib{MR3951024}{article}{
   author={S\'{a}ez, Mariel},
   author={Valdinoci, Enrico},
   title={On the evolution by fractional mean curvature},
   journal={Comm. Anal. Geom.},
   volume={27},
   date={2019},
   number={1},
   pages={211--249},
   issn={1019-8385},
   review={\MR{3951024}},
   doi={10.4310/CAG.2019.v27.n1.a6},
}

\bib{MR3090533}{article}{
   author={Savin, Ovidiu},
   author={Valdinoci, Enrico},
   title={Regularity of nonlocal minimal cones in dimension 2},
   journal={Calc. Var. Partial Differential Equations},
   volume={48},
   date={2013},
   number={1-2},
   pages={33--39},
   issn={0944-2669},
   review={\MR{3090533}},
   doi={10.1007/s00526-012-0539-7},
}

\bib{MR3133422}{article}{
   author={Savin, Ovidiu},
   author={Valdinoci, Enrico},
   title={Density estimates for a variational model driven by the Gagliardo
   norm},
   language={English, with English and French summaries},
   journal={J. Math. Pures Appl. (9)},
   volume={101},
   date={2014},
   number={1},
   pages={1--26},
   issn={0021-7824},
   review={\MR{3133422}},
   doi={10.1016/j.matpur.2013.05.001},
}
		
\end{biblist}
\end{bibdiv}
\end{document}